\documentclass{article}

\usepackage{amsmath}
\usepackage{amssymb} 
\usepackage{amsthm}  
\usepackage{amscd}
\usepackage{eucal}
 
\usepackage{fancyhdr}
\usepackage{graphicx} 
\usepackage{makeidx}  
\makeindex

\usepackage{chngcntr}
\counterwithin{figure}{section}

\usepackage{tikz}
\usetikzlibrary{calc,intersections,through,backgrounds,arrows.meta, decorations.pathmorphing, positioning}
 
  \usepackage{color}
   
\newcommand{\green}[1]{\textcolor{green}{#1}}

\newcommand{\diver}{\mathrm{div}\, }     %%%%
\newcommand{\grad}{\mathrm{grad}\, }
\newcommand{\image}{\mathrm{im\, }}
\newcommand{\simp}{simplicial complex }
\newcommand{\Vol}{\mathrm{Vol}\, }   
\newcommand{\vol}{\mathrm{Vol}\, }
\DeclareMathOperator*{\sgn}{\ensuremath{sgn}}
\DeclareMathOperator*{\conv}{\ensuremath{conv}}
\renewcommand{\vec}[1]{\mbox{\boldmath \small $#1$}}

\newcommand{\R}{{\mathbb R}}
\newcommand{\N}{{\mathbb N}}
\newcommand{\Z}{{\mathbb Z}}

\newcommand{\Cal}{\ensuremath{\mathcal{C}}}

\theoremstyle{theorem}
\newtheorem{theorem}{Theorem}[section]
\newtheorem{Thm}{Theorem}[section]

\newtheorem{lemma}{Lemma}[section]
\theoremstyle{definition}
\newtheorem{defi}{Definition}[section]

\theoremstyle{remark}

\theoremstyle{remark}

\theoremstyle{remark}
\newtheorem{remark}{Remark}[section]

\newcommand{\bel}[1]{\begin{equation}\label{#1}}

\newcommand{\be}{\begin{equation}}

\newcommand{\ba}{\begin{eqnarray}}
\newcommand{\ea}{\end{eqnarray}}
\newcommand{\rf}[1]{(\ref{#1})}
 
\newcommand{\qe}{\end{equation}}
\newcommand{\ee}{\end{equation}}

\title{From the discrete to the continuous, from simplicial complexes to Riemannian manifolds. Approximating flows and cuts on manifolds by discrete versions}
\author{Marzieh Eidi$^{2,1}$, J\"urgen Jost$^{1,2}$, Dong Zhang$^3$\\
$^1$ Max-Planck Institut for Mathematics in the Sciences, Leipzig, Germany\\%
    $^2$ ScaDS.AI, Dresden/Leipzig, Germany\\
    $^3$ School of Mathematical Sciences, Peking University, Beijing, China\\
    [2ex]%\\
        meidi@mis.mpg.de, jost@mis.mpg.de, dongzhang@math.pku.edu.cn}

\date{}

\begin{document}

\maketitle
\noindent {\bf Dedicated to Shing-Tung Yau for the celebration of his 75th birthday and in recognition of his profound and wide ranging mathematical influence }

\begin{abstract}
    Many fundamental structures of Riemannian geometry have found discrete counterparts for graphs or combinatorial ones for simplicial complexes. These include those discussed in this survey, Hodge theory, Morse theory, the spectral theory of Laplace type operators and Cheeger inequalities, and their interconnections. This raises the question of the relation between them, abstractly as structural analogies and concretely what happens when a graph constructed from random sampling of a Riemannian manifold or a simplicial complex triangulating such a manifold converge to that manifold. We survey the current state of research, highlighting some recent developments like Cheeger type inequalities for the higher dimensional geometry of simplicial complexes, Floer type constructions in the presence of periodic or homoclinic orbits of dynamical systems  or the disorientability of simplicial complexes.

\end{abstract}

\section{Introduction}
The spectral theory  of linear operators was developed by David Hilbert \cite{Hilbert04} (who also coined the name \emph{eigenvalue}, \emph{Eigenwert} in German). Specific linear operators naturally arise in geometry. After 
Bernhard Riemann had created Riemannian geometry (see \cite{Jost25}), Eugenio Beltrami \cite{Beltrami64} (see also \cite{Beltrami02}) introduced the Riemann version of the Laplace operator, now called the Laplace-Beltrami operator, that  operates on
functions on a Riemannian manifold. When the Riemannian metric tensor of the manifold in question in local coordinates $(x^1,\dots ,x^n)$ is $g_{ij}$, its inverse is $g^{ij}$, and its determinant is $g$, the Laplace-Beltrami operator operates on a function $f$ is
\begin{equation}\label{belt1}
        \Delta f = - \diver \grad f = - {1 \over \sqrt{g}} {\partial \over \partial x^j}
                \Bigl( \sqrt{g} g^{ij} {\partial f \over \partial x^i} \Bigr),
\end{equation}
where we have used the standard Einstein summation convention. 
It was first realized by Hermann Weyl \cite{Weyl12} that the eigenvalues of this operator encode important geometric information about the underlying Riemannian manifold, and this later became a very active branch of research in Riemannian geometry, with important contributions for instance by Li and Yau \cite{Li80a,Li83}. And there is the Cheeger inequality \cite{Cheeger70}, relating the first non-vanishing eigenvalue of the Laplace-Beltrami operator  to a geometric quantity that quantifies how difficult it is to cut up a Riemannian manifold into two large pieces. \\
The Laplace-Beltrami operator can be generalized to operate on differential forms. 
 When $\omega=\eta(x)dx^{i_1}\wedge \ldots \wedge dx^{i_k}$ is a $k$-form,  its differential
is the  $(k+1)$-form
\begin{equation}
  \label{h1}
d_k \eta(x)dx^{i_1}\wedge \ldots \wedge dx^{i_k}=  \sum_j \frac{\partial \eta(x)}{\partial x^j}dx^j\wedge dx^{i_1}\wedge \ldots \wedge dx^{i_k}.
\end{equation}
Using the natural Riemannian scalar product on $k$-forms, we obtain its  adjoint operator $d^\star_k$. It  maps a $(k+1)$-form to a $k$-form, and
we obtain  the Hodge Laplacian \cite{Hodge41} (which could also be named after Weyl, de Rham or Kodaira, as they all had come up about simultaneously with similar ideas)
\begin{equation}
  \label{h2}
  \Delta_k =d_{k-1}d^\star_{k-1}+ d_{k}d^\star_{k}.
\end{equation}
 For $k=0$, it is simply the Laplace-Beltrami operator \eqref{belt1}. \\
 Beno Eckmann \cite{Eckmann44} then
realized that an analogous theory can be developed for
simplicial complexes. That theory was extended in \cite{Horak13}; in particular, by considering different weight functions on the set of simplices, it was shown that there exists a one-to-one correspondence between such weight functions and the possible scalar products on the cochain group, thereby providing a unifying perspective on the various discrete Laplacians introduced in the literature.\\

Another fundamental theory is that of Morse \cite{Morse25} about the relations between the numbers and indices of the critical points of a smooth function on a Riemannian manifold and the homology of the underlying  manifold. Again, that theory has seen far reaching generalizations, in particular by Conley \cite{Conley78} who more generally related the local topology of isolated invariant sets of a dynamical systems and their exit sets with the global homology of the underlying manifold on which the dynamical system operated. Witten introduced a new perspective from the physics of supersymmetry \cite{Witten82,Witten82a} which in the hands of Floer lead to a general theory where only relative indices of critical points mattered \cite{Floer88c,Floer89}. Meanwhile,   Forman \cite{Forman98a} introduced a combinatorial version of Morse theory on cell complexes. And there are discrete versions of the Cheeger inequality for the Laplace operator of a graph. \\
These results are important intrinsically, but in the discrete case also for many applications in computer science and machine learning. They are also essentially well known. So, what can we contribute here? \\
We ask for systematic connections between these results, in particular the relation between the results in the Riemannian case and those in the discrete/combinatorial setting. These can enrich each other and lead to a more profound understanding. And this can also help to approach important open questions, like the analogues of the Cheeger type inequalities for the higher order Laplacians, those of Hodge or Eckmann type. Or we can ask about the relation between the theory of Witten and Floer and that of Conley, both in the continuous and in the discrete setting. We can report progress in these directions, but also have to point out open problems. \\
Other topics that are of great interest can be only briefly mentioned here, but not addressed in depth, like the relation between higher order Laplacians and random walks.

\section{Structures and questions}
\subsection{Simplicial complexes}\label{simp}

\subsubsection{The topology of simplicial complexes}
We consider a \simp $\Sigma$. It will be assumed finite whenever convenient,  built on a set of $N$ vertices. We denote the collection of its $k$-dimensional simplices by $\Sigma_k$. Thus, a $k$-simplex is a collection of $k+1$ vertices. As required for a simplicial complex, whenever $\sigma \in \Sigma_k$, then all its facets are contained in $\Sigma_{k-1}$. When $\sigma_k=(v_0,\dots ,v_k)$, its facets are the $(k-1)$-simplices of the form $\rho_{k-1}=(v_0,\dots ,\hat{v}_i,\dots ,v_k)$ for $i=0,\dots ,k$. We also write $\rho_{k-1}\in \partial \sigma_k$. We observe that when two different $(k-1)$-simplices $\rho_{k-1}, \rho'_{k-1} \in \partial \sigma_k$, then there can be no other $\sigma'_k$ with $\rho_{k-1}, \rho'_{k-1} \in \partial \sigma'_k$.\\
$n$, the largest $k$ for which $\Sigma_k\neq \emptyset$, is called the dimension of $\Sigma$. Obviously, $n\le N-1$. We shall implicitly assume that the \simp $\Sigma$ is connected in the sense that the graph $\Sigma_1$ is connected.\\

We shall need signs, and therefore need to consider \emph{orientations}. A $k$-simplex can carry two orientations, arbitrarily called positive and negative. They come from the ordering of its vertices. Two orientations agree (disagree) when the orderings differ by an even (odd) permutation. 
From now on, $[\sigma_k]$ will stand for an \emph{oriented} simplex, and the opposite orientation will be denoted by $-[\sigma_k]$. But in the sequel, we shall sometimes leave out the bracket $[.]$ when it should be clear from the context that we consider oriented simplices. 
The $k$-th chain group $C_k$ of $\Sigma$ is a vector space with the basis of positively oriented $k$-simplexes. The {\em boundary map } 
$\partial_d: C_k \rightarrow C_{k-1}$ is a linear operator defined by 
\begin{equation} 
\partial_k (v_0,...,v_k)=\sum_{j=0}^{k} (-1)^j (v_0,...,v_{j-1},v_{j+1},...,v_k)
\end{equation}
where we have already left out the brackets.\\
A simplex induces orientations on its facets. When 
$[\sigma_k]=(v_0,\dots ,v_k)$ is positively oriented, then so is  the facet  $[\rho_{k-1}]=(v_0,\dots ,\hat{v_i}, \dots, v_k)$. But when  $\rho_{k-1}$ is the facet of another simplex, it may or may not get the same orientation. We call $\Sigma$ \emph{orientable} (\emph{disorientable} \cite{Eidi24a}) if whenever an $(n-1)$-simplex is a facet of more than one $n$-simplex, then the induced orientations disagree (agree). A necessary, but not a sufficient condition for $\Sigma$ to be orientable, is that any $(n-1)$-simplex has degree at most 2, where the degree is the number of its cofaces, that is the number of $n$-simplices of which it is a facet.\\
For disorientability, the situation is different. This will be clarified in Section \ref{disorient} below.

\subsubsection{The Eckmann Laplacian}
We summarize the theory of \cite{Eckmann44,Horak13}, partly using the presentation of \cite{Jost/Mulas/Zhang26}. \\
We consider  the collection $C^k$  of  linear functions $\phi$ on chains of oriented simplices from $\Sigma$ that are alternating, 
\bel{lap0}
\phi(-[\sigma_k])=-\phi([\sigma_k]).
\qe
On $C^k$, the dual of $C_k$, we have the
coboundary operator $\delta_{k}$, the dual of the simplicial boundary operator $\partial_k$. It satisfies the basic relation
\begin{equation}\label{lap1}
  \delta_{k}\circ \delta_{k-1} =0.  
\end{equation}
For    (positive definite) scalar products $(\cdot,\cdot)_k$
 on the  $C^k$, we obtain 
 the adjoint
$(\delta_{k})^{*}:C^{k+1}\rightarrow ~C^{k}$  of  $\delta_{k}$  as 
 \bel{h4}
 (\delta_{k}f_{1},f_{2})_{k+1}=(f_{1},(\delta_{k})^{*} f_{2})_{k},
 \qe
for  $f_{1}\in C^{k}$ and $f_{2}\in C^{k+1}$. \eqref{lap1}
 then yields
 \begin{equation}\label{lap1a}
  {\delta_{k-1}}^*\circ {\delta_{k}}^* =0.  
\end{equation}
 In diagrammatic form
\begin{equation} \label{lap1b}
C^{k-1}
% use packages: array
\begin{array}{l}
\underrightarrow{\delta_{k-1\textrm{ }}}\\ 
 \overleftarrow{ {\delta_{k-1}}^*}
\end{array}
C^{k} \begin{array}{l}
\underrightarrow{\textrm{ }\textrm{ }\delta_{k\textrm{ }\textrm{ }}}\\ 
 \overleftarrow{ {\textrm{ }\textrm{ }\delta_{k}}^{* \textrm{ }}}
\end{array}
 C^{k+1}
\end{equation}
where the composition of two arrows in the same direction yields 0.\\
We easily observe from \eqref{h4}
  \ba
  \label{lap1.1}
  (\image \delta_k)^\perp &=& \ker  (\delta_{k})^{*},\\
  \label{lap1.2}
  (\image (\delta_k)^{*})^\perp &=& \ker  \delta_{k},
  \ea
and 
  \bel{lap2.1}
  \image \delta_{k-1}\cap  \image (\delta_k)^{*} =\{0\}, 
  \qe
  since if for instance  $f= \delta_{k-1}h=(\delta_k)^{*}g$,
  then
  \begin{equation*}
    (f,f)=(\delta_{k-1}h,(\delta_k)^{*}g)=(\delta_k \delta_{k-1}h,g)=0
  \end{equation*}
 and so  $f=0$. Thus 
  \bel{lap3.1}
  C^k=\image \delta_{k-1}\oplus  \image (\delta_k)^{*} \oplus (\ker
  \delta_{k} \cap \ker  (\delta_{k-1})^{*}).
  \qe
 Eckmann's theorem then says that
 \begin{equation}\label{eck1}
 \dim (\ker
  \delta_{k} \cap \ker  (\delta_{k-1})^{*}) =\dim H^k(\Sigma,\R)=b_k , 
 \end{equation}
 the $k$th Betti number of $\Sigma$.\\
Composing the  $\delta$ and $\delta_*$ operators yields   three self-adjoint operators, the \emph{$q$-dimensional up Laplace operator}
\begin{equation}\label{eck12}
    {L_{k}}^{up} =(\delta_{k})^{*}\delta_{k},
\end{equation}
the \emph{$q$-dimensional  down Laplace operator} 
\begin{equation}\label{eck13}
{L_{k}}^{down} =\delta_{k-1}(\delta_{k-1})^{*},
\end{equation}
and the \emph{$q$-dimensional Laplace operator}
\begin{equation}\label{eck14}
L_{k}  = {L_{k}}^{up}+ {L_{k}}^{down} =(\delta_{k})^{*}\delta_{k}+\delta_{k-1}(\delta_{k-1})^{*}.
\end{equation}
Eckmann's theorem \eqref{eck1} then says that
 \begin{equation}\label{eck2}
 \dim \ker L_k=b_k . 
 \end{equation}
Of course, for $k=0$, 
\bel{lap21}
L_0={L_{0}}^{up}.
\qe
Similarly, when $n$ is the dimension of $\Sigma$,
\bel{lap22}
L_{n}={L_{n}}^{down}.
\end{equation}
  These  operators   are non-negative, 
\bel{lap23}
(Lf,f)\ge 0 \text{ for all } f\in C^k 
\qe
for $L={L_{k}}^{up}$, ${L_{k}}^{down}$, 
and 
  \begin{eqnarray}
\label{lap27a}
{L_k}^{up}f=0 &\text{ iff }& \delta_{k}f=0\\
\label{lap27b}
{L_k}^{down}f=0 &\text{ iff }& (\delta_{k-1})^\ast f=0\\
\label{lap27}
L_kf=0 &\text{ iff }& \delta_{k}f=0\text{ and } (\delta_{k-1})^\ast f=0.
\end{eqnarray}
Since the operators ${L_{k}}^{up}$, ${L_{k}}^{down}$ and  $L_{k}$  are self-adjoint,  non-negative operators on finite-dimensional Hilbert spaces, their  eigenvalues are real and  non-negative. \\

Since ${L_{k}}^{up} =(\delta_{k})^{*}\delta_{k}$ and ${L_{k+1}}^{down} =(\delta_{k}){\delta_{k}}^{*}$, their non-zero eigenvalues are the same. In particular, if we want to get the eigenvalues of the family $L_k$ for $k=0,\dots ,n$, it suffices to determine the multiplicities of the eigenvalue 0 and the eigenvalues of one of the families, ${L_{k}}^{up}$ or ${L_{k}}^{down}$.\\

By Eckmann's theorem \eqref{eck2}, the multiplicity of the eigenvalue 0 of $L_k$ equals the Betti number $b_k$. \\
We may then naturally ask what geometric information is encoded in the other eigenvalues, those that do not vanish.

\subsection{Riemannian manifolds}
The theory for Riemannian manifolds is completely analogous to that for simplicial complexes. We consider Riemannian manifolds that, for simplicity, we assume to be connected (to make the eigenvalue 0 of the Laplace-Beltrami operator simple) and compact (to avoid certain analytical difficulties that are not central for our story).  \\
Analogously to  \eqref{lap3.1}, we have the  Hodge decomposition for the Hodge Laplacian in Riemannian geometry. From the exterior differential $d_k$ and its adjoint $d_k^*$, which satisfy
\begin{equation}
  d_k \circ d_{k-1}=0 =  d_{k-1}^*\circ d_k^*, 
\end{equation}
we can form the three operators, the \emph{$k$-dimensional up Laplace operator}
\begin{equation}\label{lap11}
   {\Delta_{k}}^{up} =(d_{k})^{*}d_{k}, 
\end{equation}
the \emph{$k$-dimensional  down Laplace operator}
\begin{equation}\label{lap12}
 {\Delta_{k}}^{down} =d_{k-1}(d_{k-1})^{*},   
\end{equation}
and the \emph{$k$-dimensional Laplace operator} \eqref{h2}, 
\begin{equation}\label{lap13}
   \Delta_{k}  = {\Delta_{k}}^{up}+ {\Delta_{k}}^{down} =(d_{k})^{*}d_{k}+d_{k-1}(d_{k-1})^{*}. 
\end{equation}
Since the operators ${\Delta_{k}}^{up}$, ${\Delta_{k}}^{down}$ and  $\Delta_{k}$  are self-adjoint,  non-negative operators on  Hilbert spaces, their  eigenvalues are real and  non-negative, and when the underlying Riemannian manifold is compact, as we assume here, the spectrum is discrete. \\

Again, the non-zero eigenvalues of ${\Delta_{k}}^{up}$ and ${\Delta_{k+1}}^{down} $ are the same. So, if we want to get the eigenvalues of the family $\Delta_k$ for $k=0,\dots ,n$, where $n$ is the dimension of the Riemannian manifold in question, it suffices to determine the multiplicities of the eigenvalue 0 and the eigenvalues of one of the families, ${\Delta_{k}}^{up}$ or ${\Delta_{k}}^{down}$.\\

By the Hodge theorem \eqref{eck2}, the multiplicity of the eigenvalue 0 of $L_k$ equals the Betti number $b_k$,
\begin{equation}\label{h5}
    \dim \ker \Delta_k =b_k.
\end{equation}
Again, we ask what geometric information is encoded in the non-zero eigenvalues.

\subsection{Rayleigh quotients}\label{rayleigh}
As discovered by Courant-Fischer and Weyl, the eigenvalues of all the symmetric operators $A$ considered here (Eckmann Laplacians on finite simplicial complexes, Laplace-Beltrami and Hodge Laplacians on compact Riemannian manifolds) can be obtained as the critical values of the Rayleigh quotients
\begin{equation}\label{r1}
    \frac{(Af,f)}{(f,f)}
\end{equation}
for functions $f$ with $(f,f)\neq 0$. For our Laplace operators $A=\delta \circ \delta^* +\delta^* \circ \delta$, the scalar product is of course the one used for the definition of the adjoint $\delta^*$, and \eqref{r1} becomes
\begin{equation}\label{r2}
    \frac{(\delta f, \delta f)+(\delta^* f,\delta^* f)}{(f,f)}.
\end{equation}
The critical points are the eigenfunctions, and eigenfunctions $u_i,u_j$ for eigenvalues $\lambda_i,\lambda_j$  are orthogonal to each other, $(u_i,u_j)=0$, and in fact, we can find a basis of eigenfunctions satisfying
\begin{equation}\label{r3}
  (u_i,u_j)=\delta_{i,j}.  
\end{equation}
For the Laplace-Beltrami operator, and analogously the Eckmann Laplacian $L_0$ on vertices, the smallest eigenvalue, that is, the smallest critical value of \eqref{r1}, is $\lambda_1=0$, with a constant eigenfunction $u_1=c$. Therefore, the higher eigenfunctions satisfy
\begin{equation}\label{r4}
  (u_j,1)=0\quad \text{ for } j=2,\dots  .
\end{equation}
In the manifold case, this means 
\begin{equation}\label{r5}
\int_M u_j =0\quad  \text{ for } j=2,\dots  ,
    \end{equation}
while in the discrete case below, we shall have
\begin{equation}\label{r6}
\sum_v \deg v\ u_j(v) =0\quad  \text{ for } j=2,\dots  .
    \end{equation}

\subsection{$p$-Laplacians}\label{plap} 

Below, we shall need a non-linear generalization of Laplacians, the $p$-Laplacians that describe nonlinear diffusions whose strengths depend on the magnitude of the gradient. The Laplacians that we have considered so far represent the case $p=2$, because we have a square in \eqref{r2}. When we take the numerator to the power $p/2$, we obtain the $p$-Laplacian. The case $p=1$ presents some subtleties, but this will be precisely the case that we need, because the Cheeger constant (that will be introduced later) is an
$L^1$-quantity.

Thus, for $p>1$,
we put $$a_p:(t_1,t_2,\cdots)\mapsto
(|t_1|^{p-2}t_1,|t_2|^{p-2}t_2,\cdots).$$
This becomes undetermined for
$p=1$ when $t=0$, as $\frac{t}{|t|}$ doesn't possess a definite limit for $t\to 0$. Therefore, we let it be set valued for $t=0$, to include the entire set $[-1,1]$
that is, $$a_1:(t_1,t_2,\cdots)\mapsto
\left\{(\xi_1,\xi_2,\cdots):\xi_i\in\mathrm{Sgn}(t_i):=\begin{cases}
 \{1\} & \text{if } t_i>0,\\
 [-1,1] & \text{if }t_i=0,\\
 \{-1\} & \text{if }t_i<0
 \end{cases}\quad \right\}.$$
We  define 
 the $d$-th up $p$-Laplacian
 $$L^{up}_{d,p}:=\delta_d^*a_{p}\delta_d,$$
 and similarly 
 the $d$-th down $p$-Laplacian
 $L^{down}_{d,p}:=\delta_{d-1}a_p\delta_{d-1}^*$
 and  the $d$-th
  $p$-Laplacian as
 $ L_{d,p}:=L^{up}_{d,p}+ L^{down}_{d,p}$.
We can then formulate the eigenvalue problem for $L^{up}_{d,p}$ as finding $\lambda$ and nonzero functions $f:\Sigma_d\to\R$ with 
\[L^{up}_{d,p}f=\lambda a_p(f) \text{ if }p>1,\]
or  
\[0\in L^{up}_{d,1}f-\lambda a_1(f) \text{ if }p=1.\]
For $p\neq 2$, these are non-linear problems. \\
We will be interested in those eigenvalues that
can be obtained from Rayleigh quotients. 
We note that for  finite graphs, for the non-linear problem ($p\neq 2$), and in particular for $p=1$, in contrast to the linear problem $p=2$, the number of eigenvalues can  exceed the dimension of the space $C^d$ \cite{Amghibech03,Deidda23}. Classical results in nonlinear analysis guarantee that one can characterize a set of eigenvalues, whose number, counted with multiplicity, is equal to the dimension of $C^d$. These are the so-called minimax eigenvalues:
\begin{equation}\label{eq:discrete}
\lambda_k(L^{up}_{d,p}):=\inf_{\gamma(S)\ge k}\sup_{f\in S}\frac{\|\delta_d f\|^p}{\|f\|^p} ,\;\; k=1,2,\ldots
\end{equation}
where $\gamma(S)$ is the \emph{Krasnoselskii genus} of an origin-symmetric compact set $S$ defined by $$\gamma(S):=\min\{k\in\mathbb{Z}^+:\exists\text{ odd continuous }\varphi:S\to \mathbb{R}^k\setminus\{0\}\}.$$
It is known that $\lambda_2(L^{up}_{0,1})$ is the first nontrivial eigenvalue of the graph 1-Laplacian $L^{up}_{0,1}$, which equals the Cheeger constant on graphs. There is a monotonicity theorem that refines and generalizes the Cheeger inequality \cite{Zhang25}. 
\begin{theorem}
Suppose $d=0$. For any $k$, the $k$-th minimax eigenvalue $\lambda_k(L^{up}_{0,p})$ is locally Lipschitz continuous with respect to $p$, and moreover, 
\begin{itemize}
    \item the function $p\mapsto p(\frac{2}{D}\lambda_k(L^{up}_{0,p}))^{\frac1p}$ 
is increasing  on $%p\in
[1,+\infty)$, where $D:=\max\limits_{v\in\Sigma_0}\deg v$ and $\deg v$ denotes the number of $1$-simplexes containing $v$. 
\item %and 
the function 
$p\mapsto2^{-p}\lambda_k(L^{up}_{0,p})$ is decreasing on $[1,+\infty)$.
\end{itemize}  
\end{theorem}
This theorem solves an open problem posed by Amghibech \cite[Question 2]{Amghibech03}.

\subsection{Lov\'asz extension}\label{lovasz}
We will use a tool for systematically converting discrete into continuous quantities, the Lov\'asz extension. 
We consider a finite set $V=\{1,\cdots,n\}$ and its power set $\mathcal{P}(V)$, and a function 
$F:\mathcal{P}(V)\to \R$, usually assuming $F(\varnothing)=0$.\\
We identify $A\in \mathcal{P}(V){\setminus\{\varnothing\}}$ with its indicator
vector $\vec1_A\in \R^V=\R^n$  ($(\vec1_A,e_j)=1$ iff $j\in A$, where
  $e_j:=j$th unit vector). The \emph{Lov\'asz extension}  extends the domain of $F$ to the whole Euclidean space. Precisely, for any $\vec x\in\R^V$,
 \begin{align}\label{eq:Lovaintegral}
F^{L}(\vec x):&=\int_{\min\limits_{1\le i\le n}x_i}^{\max\limits_{1\le i\le n}x_i} F(V^t(\vec x))d t+F(V)\min_{1\le i\le n}x_i \\
             &= \int_{-\infty}^0(F(V^t(\vec x))-F(V))dt+ \int_0^{+\infty}(F(V^t(\vec x)){-F(\varnothing)}) d t,
 \end{align}
%  and if we add the natural assumption $F(\varnothing)=0$,
% \begin{equation}\label{eq:Lovaintegral2}f^{L}(\vec x)=\int_{-\infty}^0(F(V^t(\vec x))-F(V))dt+ \int_0^{+\infty}F(V^t(\vec x)) d t,
% \end{equation}
where $V^t(\vec x)=\{i\in V:  x_i>t\}$.

$F:\mathcal{P}(V)\to \R$   is \emph{submodular} if $F(A)+F(B)\ge F(A\cup B)+F(A\cap B)$, $\forall A,B\in\mathcal{P}(V)$. The Lov\'asz extension turns a submodular function into a convex function, and we can hence minimize the former by minimizing the latter:
\begin{theorem}[Lov\'asz]
	$F:\mathcal{P}(V)\to\mathbb{R}$ is submodular iff  $F^L$ is convex.
\end{theorem}

\begin{theorem}[Lov\'asz]If $F:\mathcal{P}(V)\to\mathbb{R}$ is submodular with $F(\varnothing)=0$, then
	$$\min\limits_{A\subset V}F(A)=\min\limits_{\vec x\in [0,1]^V}F^L(\vec x).$$
\end{theorem}

As in \cite{Jost/Zhang21a,Jost/Zhang21b}, we  now introduce the \emph{$(F_1,F_2)$--constant} corresponding to a given pair of functions $F_1:\mathcal{P}(V)\to\R$  and $F_2:\mathcal{P}(V)\to[0,+\infty)$ with $F_2(A)>0$ whenever $A\ne\emptyset,V$, which is defined as \begin{equation}\label{eq:Cheeger}
h(F_1,F_2):= \min\limits_{A\in\mathcal{P}(V)\setminus\{\emptyset,V\}}\frac{F_1(A)}{\min\{F_2(A),F_2(V\setminus A)\}} .   \end{equation}
This will be used in our discussion of the Cheeger constant below, and we have \cite{Jost/Mulas/Zhang26}
\begin{theorem}[]
$$h(F_1,F_2)=\min\limits_{x\in\R^n: x\text{ non-constant}}\,\frac{F_1^L(x)}{\min\limits_{t\in\R} F_2^L(| x-t\vec 1|)}$$
where $| x-t\vec 1|:=(|x^1-t|,\ldots,|x^n-t|)\in\R^n$.
\end{theorem}

\subsection{Analogies and convergence}
In the preceding sections, we have described the Hodge theory for Laplacians on Riemannian manifolds and the Eckmann theory for Laplacians on simplicial complex. We have seen an obvious structural analogy. But is there more? \\

Riemannian manifolds can be approximated by simplicial complexes in various ways. We can triangulate a Riemannian manifold. The simplices then carry a geometric structure induced from the Riemannian manifold. Eckmann's theory, however, depends only on the combinatorial structure and a choice of scalar products for functions on the $k$-dimensional chains. Thus, when we make the triangulation finer and finer, how much will the combinatorial structure reveal about the metric one, and what is the best choice of scalar products?
 There is work by Dodziuk and Patodi in this direction \cite{Dodziuk76,DodziukPatodi76} which we shall take up below from a modern perspective.\\
 
 We can also approximate a manifold by a graph. That can be done, for instance, by randomly sampling points from the manifold, according to its volume measure, and connecting points that are sufficiently close. We can then ask whether the Laplacian of these graphs and their eigenvalues converge to those of the Riemannian manifold, perhaps after some rescaling. This is the basis of the method of Belkin and Niyogi \cite{Belkin03,Belkin08} that has become one of the most important techniques in machine learning (for the mathematical analysis, see for instance also \cite{Joharinad23}). But this applies then only for the $0$-Laplacians, that is, for the Laplace-Beltrami operator of the manifold and the graph Laplacian. Can this be extended to higher order Laplacians?
 \begin{remark} 
 The normalized Laplacian of a graph
 \begin{equation}
     \Delta f(u)=f(u)-\sum_{v\sim u}\frac{w_{uv}}{\deg u}f(v)
 \end{equation}
 is the infinitesimal generator of a random walk on the graph.
 It encodes the probabilities $\frac{w_{uv}}{\deg u}$ with which a walker at $u$ will reach a vertex $v$ in the next step. Conversely (as explained for instance in \cite{Jost14}), the operator
 \begin{equation}
     \Delta' g(v)=g(v)-\sum_{u\sim v}\frac{w_{uv}}{\deg u}g(u)
 \end{equation}
 encodes the probability from which vertex $u$ a walk landing at $v$ arrived. These two operators have the same eigenvalues, but an eigenfunction $h(v)$ of $\Delta$ is changed into the eigenfunction $\deg v\ h(v)$ of $\Delta'$.\\
 As mentioned above, $\Delta$ serves as a discrete approximation to the (weighted) Laplace–Beltrami operator on manifolds, which itself is the infinitesimal generator of Brownian motion \cite{Hsu2002}. \\
 For more details on this convergence and its rate, see, for instance, \cite{Hein2007}. (Lazy) random walks on connected graphs are irreducible and their limiting behavior is independent of the starting vertex \cite{Lovasz1993}; however, such irreducibility is generally absent in Laplacian-based random walks on higher-dimensional simplicial complexes, as analyzed and demonstrated in \cite{Rosenthal17, Mukherjee2016, Eidi23}. This structural difference is key for understanding their limiting and convergence behavior. Moreover, there is no fully developed stochastic process whose generator coincides with the Hodge Laplacian on general 
$k$-forms on manifolds. We omit the details here, but there remain many open problems in this direction for interested researchers.   
\end{remark}

\section{Morse theory and flows}
\subsection{Morse theory on Riemannian manifolds}
Morse theory tells us that the number $m_k$ of critical points of index $k$ of a smooth function $f$ on a compact Riemannian manifold $M$ ($\dim M=n$, metric tensor $g_{ij}$) that has only non-degenerate critical points satisfies
\begin{equation}\label{m1}
m_k\ge b_k
    \end{equation}
    and more generally
\begin{equation}\label{m1a}
m_k-m_{k-1}\pm \dots (-1)^k m_0 \ge b_k-b_{k-1}\pm \dots (-1)^k b_0
    \end{equation} 
    and 
    \begin{equation}\label{m2}
\sum_{i=0}^n (-1)^i m_k=    \sum_{i=0}^n (-1)^i b_k ,   
    \end{equation}
where the index of a critical point is the number of negative eigenvalues of the Hessian of $f$ at that point and $b_k$ is the $k$-th Betti number. 
\\
More generally, we have for $t\in \R$
\begin{equation}\label{m3}
  \sum_{i=0}^n m_k t^k   - \sum_{i=0}^n  b_kt^k =(1+t)   \sum_{i=0}^n q_k t^k  
\end{equation}
where the $q_k \in \N$. $t=-1$ in \eqref{m3} yields \eqref{m2}.

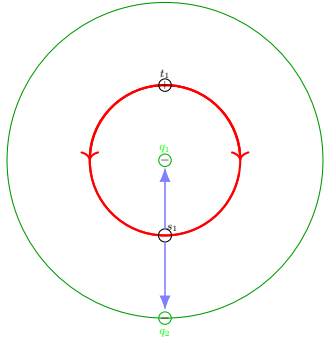
\begin{figure}
\centering
\begin{tikzpicture}[scale=0.5, transform shape,
    every node/.style={font=\small},
    flow/.style={-Latex, thick},
    addflow/.style={-Latex, blue!50, line width=0.6pt},
    saddle/.style={circle, draw=black, inner sep=3.5pt},
    sink/.style={circle, draw=green!70!black, inner sep=0pt, minimum size=7pt},
    top/.style={circle, draw=black!70!black, inner sep=0pt, minimum size=7pt},
    closed/.style={red, thick},
  arrowcircle/.style={
    postaction={decorate},
    decoration={markings, mark=at position 0.25 with {\arrow{>}}}
  }
]

%--- outer disk (S^2 as disk model) ---
\draw[green!60!black, line width=.4pt] (0,0) circle (4.2cm);

%--- red closed orbit (γ) near top-left quadrant ---
\draw[closed] (0,0) circle (2cm);

\draw[<-,red, thick] (2,0) arc (0:350:2cm);

\draw[<-,red,thick] (-2,0) arc[start angle=180, end angle=0, radius=2cm];

%--- top (plus signs) ---
\node[top] (t1) at (0,2) {\(+\)};
\node[black] at (0,2.3) {$t_1$};
%--- saddles ---
\node[saddle] (s1) at (0,-2) {};
\node[black] at (0.2,-1.8) {$s_1$};
%--- sinks (minus signs) ---
\node[sink] (q1) at (0,0) {\(-\)};
\node[green] at (0,0.3) {$q_1$};
\node[sink,label=below:\green{$q_2$}] (q2) at (0,-4.2) {\(-\)};

%--- additional faint flow lines (blue) for visual context ---
\draw[addflow] (s1) .. controls (0,-3)  .. (0,-4);

\draw[addflow] (s1) .. controls (0,-1)  .. (q1);
\end{tikzpicture}
\caption{A Morse function on $S^2$}\label{fig1}
\end{figure}

Fig. \ref{fig1} indicates a Morse function on $S^2$. This is of course a trivial situation, but we shall use it to explain some points, and as a reference for future figures.\\
The green circle represents a single point $q_2$. The function has one maximum at $t_1$ (index 2), a saddle at $s_1$ (index 1) and two minima at $q_1,q_2$ (index 0). Thus, $m_2-m_1+m_0=2$, the Euler number of $S^2$. We shall now count flow lines modulo 2, although we could also do an integer count, if we take orientations into account. As we shall explain below, counting flow lines between critical points of index difference 1 yields a boundary operator. There are two flow lines from $t_1$ to $s_1$, thus mod 2, we have $\partial_2 t_1 =0$. There is one flow each from $s_1$ to $q_1$ and $q_2$,
hence $\partial_1 s_1=q_1+q_2$. And since the latter are minima, no flow lines emanate from them, hence $\partial_0 q_1=0=\partial_0 q_2$. Therefore, at index 2, $\dim \ker \partial_2=1$; at index 1 $\dim \ker \partial_1 =0$ and $\dim \image \partial_2 =0 $; and at index 0, $\dim \ker \partial_0 =2$ and $\dim \image \partial_1 =1$. Thus, $b_2=1,b_1=0, b_0=1$, the Betti numbers of $S^2$.

\subsection{The Witten Laplacian}\label{witten}
Witten's idea \cite{Witten82}  was to recover these inequalities
from the analytic behavior of a deformed Laplacian
and to interpret Morse theory as a semiclassical limit
of a supersymmetric quantum system. He noticed that the de~Rham complex \( (\Omega^*(M), d) \)
can be \emph{deformed} by the function \(f\).
For a real parameter \( t\ge 0 \), define
\begin{equation}\label{eq:dt}
d_t = e^{-tf}\, d\, e^{tf},
\qquad
d_t^* = e^{tf}\, d^*\, e^{-tf}\qquad \text{ for }t\ge0.
\end{equation}
This is a conjugation of \(d\) by the weight \(e^{tf}\),
and hence \(d_t\) is cohomologically equivalent to \(d\).
Intuitively, this means that we are ``tilting'' the geometry by the potential~\(f\):
the forms are now measured in a weighted space that emphasizes regions
where \(f\) is small.
With the modified differential \(d_t\), one defines the \emph{Witten Laplacian}
\begin{equation}\label{eq:Ht}
H_t = d_t^*d_t + d_t d_t^* 
\end{equation}
and  of course can do this also for the up- and down-Laplacians.\\
The operator $H_t$ acts on differential forms of all degrees.
A straightforward calculation shows that \(H_t\) can be expanded as
\begin{equation}\label{eq:Ht-expanded}
H_t = \Delta + t^2|\nabla f|^2 + t\, \mathrm{Hess}(f)\text{-terms},
\end{equation}
where \( \Delta = dd^* + d^*d \) is the Hodge Laplacian
and the last term involves contractions of the Hessian of \(f\)
with the form components.
Thus \(H_t\) is a Schrödinger-type operator:
the Laplacian plus a potential term \(t^2|\nabla f|^2\).
In the physics analogy, it plays the role of a \emph{Hamiltonian}.
When \(t\) becomes large, the potential \(t^2|\nabla f|^2\) dominates everywhere
except near the critical points of \(f\) (where \(\nabla f=0\)).
Consequently, the lowest eigenstates of \(H_t\), those with very small eigenvalues,
become strongly \emph{localized} around the critical points of \(f\).
In the limit \(t\to\infty\),
each critical point of index \(k\) contributes roughly one localized eigenstate
of degree \(k\).
Therefore, the number of small eigenvalues of \(H_t\) acting on \(k\)-forms
is approximately \(m_k\) from which Witten obtains the Morse inequalities \eqref{m1}, \eqref{m1a}, \eqref{m2}.

%Witten \cite{Witten82} then conjugated the exterior derivative (and consequently also its adjoint)
%\begin{equation}\label{m4}
   % d_{t}=e^{-tf} d e^{tf}, \quad d_t^*=e^{tf} d^* e^{-tf} \qquad \text{ for } t\ge 0 
%\end{equation}
%(where we leave out the index $k$) 
%to modify the Laplacian \eqref{lap13} as
%\begin{equation}\label{m5}
 % H_t=(d_{t})^{*}d_{t}+d_{t}(d_{t})^{*}, 
%\end{equation}
%and correspondingly of course also for the up- and down-Laplacians.\\
%For $t\to \infty$, the eigenstates for the lowest eigenvalues (essentially 0) concentrate at the critical points of $f$, from which Witten obtains the Morse inequalities \eqref{m1}.\\

Between critical points $p,q$ with Morse indices $i(q)=i(p)-1$, we can look at the gradient flow lines $c:[0,1]\to M$ with $c(0)=p, c(1)=q$ satisfying
\begin{equation}\label{m6}
  \frac{dc(t)}{dt}=\grad f(c(t)),   
\end{equation}
or in local coordinates
\begin{equation}\label{m7}
  \frac{dc^i(t)}{dt}=g^{ij}\frac{\partial f}{\partial x^j}(c(t)).  
\end{equation}
Counting these flow lines $c_\alpha, \alpha=1,\dots, m$ from $p$ to $q_\alpha$ with appropriate orientations $o(c_\alpha)=\pm 1$ yields a boundary operator 
\begin{equation}\label{m8}
\partial p= \sum_\alpha o(c_\alpha) q_\alpha
    \end{equation}
    that satisfies
    \begin{equation}\label{m9}
  \partial \circ \partial =0      
    \end{equation}
which recovers the cohomology of $M$ and is the key operator of Floer's theory \cite{Floer89}. (A systematic presentation of Floer theory is given in \cite{Schwarz93}. See also \cite{Banyaga04,Jost17}.)

\begin{figure}
\centering
\begin{tikzpicture}[scale=0.5, transform shape, 
    every node/.style={font=\small},
    flow/.style={-Latex, thick},
    addflow/.style={-Latex, blue!50, line width=0.6pt},
    saddle/.style={circle, draw=black, inner sep=3.5pt},
    sink/.style={circle, draw=green!70!black, inner sep=0pt, minimum size=7pt},
    top/.style={circle, draw=black!70!black, inner sep=0pt, minimum size=7pt},
    closed/.style={red, thick},
  arrowcircle/.style={
    postaction={decorate},
    decoration={markings, mark=at position 0.25 with {\arrow{>}}}
  }
]

%--- outer disk (S^2 as disk model) ---
\draw[green!60!black, line width=.4pt] (0,0) circle (4.2cm);

%--- red closed orbit (γ) near top-left quadrant ---
\draw[closed] (0,0) circle (2cm);

\draw[<-,red, thick] (2,0) arc (0:350:2cm);

\draw[<-,red,thick] (-2,0) arc[start angle=180, end angle=0, radius=2cm];

%--- top (plus signs) ---
\node[top,label=above:$t_1$] (t1) at (0,2) {\(+\)};

%--- saddles ---
\node[saddle] (s1) at (0,-2) {};
\node[black] at (0.2,-1.8) {$s_1$};
\node[saddle] (s2) at (0,0) {};
\node[black] at (0.2,0.2) {$s_2$};

%--- sinks (minus signs) ---
\node[sink,label=above:$q_1$] (q1) at (-1,0) {\(-\)};
\node[sink,label=above:$q_2$] (q2) at (1,0) {\(-\)};
\node[sink,label=below:\green{$q_3$}] (q3) at (0,-4.2) {\(-\)};

%--- additional faint flow lines (blue) for visual context ---
\draw[addflow] (s1) .. controls (0,-3)  .. (0,-4);
\draw[addflow] (s1) .. controls (-.2,-1)  .. (q1);
\draw[addflow] (s1) .. controls (.2,-1)  .. (q2);

\draw[addflow] (t1) .. controls (0,1)  .. (s2);
\draw[addflow] (s2) .. controls (-.5,0)  .. (q1);
\draw[addflow] (s2) .. controls (.5,0)  .. (q2);
\end{tikzpicture}
   \label{fig1a}
    \caption{An example that does not satisfy the  condition for the boundary operator} 
   \end{figure}
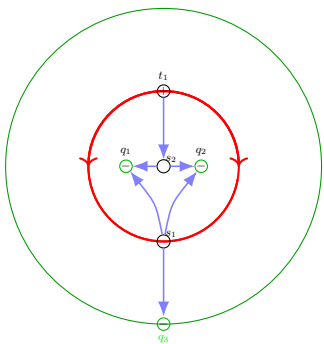
\begin{figure}
\centering
\begin{tikzpicture}[scale=0.5, transform shape,
    every node/.style={font=\small},
    flow/.style={-Latex, thick},
    addflow/.style={-Latex, blue!50, line width=0.6pt},
     thickflow/.style={-Latex, red!50, line width=1pt},
    saddle/.style={circle, draw=black, inner sep=3.5pt},
    sink/.style={circle, draw=green!70!black, inner sep=0pt, minimum size=7pt},
    top/.style={circle, draw=black!70!black, inner sep=0pt, minimum size=7pt},
    closed/.style={red, thick},
  arrowcircle/.style={
    postaction={decorate},
    decoration={markings, mark=at position 0.25 with {\arrow{>}}}
  }
]

%--- outer disk (S^2 as disk model) ---
\draw[green!60!black, line width=.4pt] (0,0) circle (4.2cm);

%--- red closed orbit (γ) near top-left quadrant ---
% \draw[closed] (0,0) circle (2cm);

% \draw[->,red, thick] (-.67,-1.9) arc (250:90:2cm);

% \draw[<-,red,thick] (2,0) arc[start angle=0, end angle=30, radius=2cm];

%--- top (plus signs) ---
\node[top,label=above:$t_1$] (t1) at (0,2) {\(+\)};
\node[top] (t2) at (0,-2) {\(+\)};
\node[black] at (0.2,-1.8) {$t_2$};

%--- saddles ---
\node[saddle,label=left:$s_1$] (s1) at (-2,0) {};
\node[saddle,label=right:$s_2$] (s2) at (2,0) {};

\node[saddle] (s3) at (0,0) {};
\node[black] at (0.2,0.2) {$s_3$};

%--- sinks (minus signs) ---
\node[sink,label=above:$q_1$] (q1) at (-1,0) {\(-\)};
\node[sink,label=above:$q_2$] (q2) at (1,0) {\(-\)};
\node[sink,label=below:\green{$q_3$}] (q3) at (0,-4.2) {\(-\)};

%--- additional faint flow lines (blue) for visual context ---
%\draw[addflow] (t2) .. controls (0,-3)  .. (0,-4);
\draw[addflow] (s1) .. controls (-1.5,0)  .. (q1);
\draw[addflow] (s2) .. controls (1.5,0)  .. (q2);
\draw[addflow] (s1) .. controls (-2.5,-1)  .. (q3);
\draw[addflow] (s2) .. controls (2.5,-1)  .. (q3);
\draw[addflow] (t1) .. controls (0,1)  .. (s3);
\draw[addflow] (t2) .. controls (0,-1)  .. (s3);
\draw[addflow] (s3) .. controls (-.5,0)  .. (q1);
\draw[addflow] (s3) .. controls (.5,0)  .. (q2);

\draw[thickflow] (t1) .. controls (-1.5,1.5) and (-1.8,1) .. (s1);
\draw[thickflow] (t1) .. controls (1.5,1.5) and (1.8,1) .. (s2);
\draw[thickflow] (t2) .. controls (-1.5,-1.5) and (-1.8,-1) .. (s1);
\draw[thickflow] (t2) .. controls (1.5,-1.5) and (1.8,-1) .. (s2);
\end{tikzpicture}
   \caption{Resolving the problem of Fig.\ref{fig1a}
      by adding a pair of critical points, a maximum and a saddle}\label{fig1b}
\end{figure}
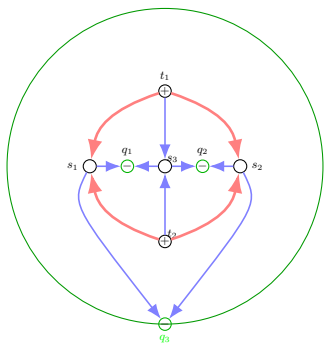

In Fig. \ref{fig1a}, we see an example where the counting of flow lines for constructing a boundary operator does not work as it is a degenerate case where there are three flow lines from $s_1$ to $q_1, q_2, q_3$ and therefore there would necessarily also be a flow line from the saddle $s_1$ to the saddle $s_2$, that is, between critical points of the same index, but that is not permitted. But this can be remedied by introducing another maximum and another saddle as in Fig. \ref{fig1b}.

\begin{remark}
Physical interpretation: Witten \cite{Witten82a} came to this construction from 
\emph{supersymmetric quantum mechanics}.
In this correspondence:
\begin{itemize}
  \item \(d_t\) and \(d_t^*\) play the role of \emph{supercharges};
  \item \(H_t\) is the \emph{supersymmetric Hamiltonian};
  \item zero-energy states (\(\ker H_t\)) correspond to
        \emph{supersymmetric ground states}, i.e., harmonic forms representing
        cohomology classes;
  \item excited states correspond to gradient flow lines between critical points.
\end{itemize}
In this way, Morse theory arises as the \emph{semiclassical limit}
of a supersymmetric system, where topology
(the Betti numbers \(b_k\)) manifests itself as the number of supersymmetric ground states. 
\end{remark}

\subsection{Morse-Smale systems and the Floer boundary operator}

From that perspective, what is fundamental are not the critical points but the boundary operator $\partial$  that counts flow lines. So, naturally, this can also be applied to flows that do not necessarily come from the gradient of a function, but from a more general dynamical system, like one that has both non-degenerate fixed points and non-degenerate periodic orbits, that is, a   Morse-Smale system. 
For such extension, as shown in \cite{Eidi21, Eidi24}, the key idea is to locally perturb periodic orbits into heteroclinic connections between two fixed points, following a result of Franks \cite{Franks79}. We show that the Floer boundary operator arising from this Franks-type replacement carries additional structure inherited from the orbits that have been perturbed away. This structure enables a systematic arrangement of flow lines between the critical points of the perturbed system, leading naturally to the definition of a boundary operator that remains meaningful even in the presence of those orbits. In other words, one can determine the boundary operator directly from the relations between the original orbits and the critical points, without explicitly invoking the perturbation, although the perturbation provides the geometric intuition for why this operator squares to zero. \\
We can not only allow for periodic orbits, but also for homoclinic ones. Suitable nondegeneracy conditions need to be satisfied, like hyperbolicity of critical points and orbits, transversal intersections between stable and unstable manifolds, and the exclusion of connecting orbits between objects of the same index, to avoid cases as in Figure \ref{fig1a} or \ref{fig4} (or perturb such cases away, as indicated in Fig. \ref{fig1b}). \\
A periodic orbit $O_{k}$ of index
$k$ generates an element $O_k^1$ of dimension $k+1$ and an element $O_{k}^0$
in dimension $k$, as for instance $S^1$ which when considered as a  closed orbit  would correspond to one object in dimension 1 and another one in dimension 0, as would be obtained when replacing it by one critical point of index 1 and another of index 0, with two flow lines from the former to the latter. Analogously for homoclinics $H_k$. Thus, each such
orbit carries topology in two adjacent dimensions and therefore corresponds to two
elements in the boundary calculus.
The differential $\partial_k : C_k(X) \longrightarrow C_{k-1}(X)$  then  counts the
number $\alpha$ of connected components of $M(\beta_i, \beta_ j )$ (mod 2) where
$\beta_i$ and $\beta_j$ are isolated rest points $p_k$ or closed homoclinic or periodic orbits.  More precisely, Thus, our boundary
operator is:
          \begin{eqnarray*}
\partial p_k &=& \sum \alpha (p_k, p_{k-1})  p_{k-1} + \sum  \alpha (p_k,  O_{k-2})   O^1_{k-2}\\ &+& \sum  \alpha (p_k,  H_ {k-2})   H_ {k-2}^1  \\\\ \partial  O^0_{k}&=& \sum \alpha ( O_{k},  O_{k-1})  O^0_{k-1}+  \sum \alpha ( O_{k},  H_ {k-1})  H_ {k-1}^0 \\ &+&  \sum \alpha ( O_{k}, p_{k-1}) p_{k-1} \\\\  \partial  O^1_{k-1}&=& \sum \alpha( O_{k-1},  O_{k-2})  O^1_{k-2} + \sum \alpha ( O_{k-1},  H_ {k-2})  H_ {k-2}^1 \\\\  \partial  H_ {k}^0&=& \sum \alpha ( H_ {k},  H_ {k-1}) H_ {k-1}^0+  \sum \alpha ( H_ {k},  O_ {k-1})  O_ {k-1}^0 \\ &+&   \sum \alpha ( H_ {k}, p_{k-1}) p_{k-1} \\\\  \partial  H_ {k-1}^1&=& \sum \alpha ( H_ {k-1}, H_ {k-2})  H_ {k-2}^1+ \sum \alpha ( H_ {k-1},  O_{k-2})  O^1_{k-2}.  
\end{eqnarray*}
Thus, not only critical points, but also orbits can occur as boundary contributions, even of critical points. 
In this definition, the sums extend over all the elements on the right hand
side; for instance, the first sum in the first line is over all critical
points $p_{k-1}$ of index $k-1$.

\begin{figure}
\centering
  \begin{tikzpicture}[scale=0.5, transform shape,
    every node/.style={font=\small},
    flow/.style={-Latex, thick},
    addflow/.style={-Latex, blue!50, line width=0.6pt},
    saddle/.style={circle, draw=black, inner sep=3.5pt},
    sink/.style={circle, draw=red!70!black, inner sep=0pt, minimum size=7pt},
    closed/.style={red, thick},
  arrowcircle/.style={
    postaction={decorate},
    decoration={markings, mark=at position 0.25 with {\arrow{>}}}
  }
]

%--- outer disk (S^2 as disk model) ---
\draw[green!60!black, line width=0.5pt] (0,0) circle (4.2cm);

%--- red closed orbit (γ) near top-left quadrant ---
\draw[closed] (0,0) circle (2cm);
\node[red] at (2,1) {$O$};
\draw[<-,red, thick] (2,0) arc (0:350:2cm);

\draw[<-,red,thick] (0,-2) arc[start angle=270, end angle=360, radius=2cm];

%--- sinks (minus signs) ---
\node[sink,label=below:$q_1$] (q1) at (0,0) {\(-\)};
\node[sink,label=above:\green{$q_2$}] (q2) at (0,4.2) {\(-\)};

%--- additional faint flow lines (blue) for visual context ---
\draw[addflow] (-2.2,0) .. controls (-2.2,0.4) and (-2.0,0.8) .. (-3.2,2);
\draw[addflow] (-.1,-2.1) .. controls (-.7,-3) and (-.4,-3.5) .. (0,-4);
\draw[addflow] (.1,2.1) .. controls (.7,3) and (.4,3.5) .. (0,4);
\draw[addflow] (2.2,-0) .. controls (2.2,-.4) and (2,-.8) .. (3.2,-2);
\draw[addflow] (-1.7,-0.8) .. controls (-1.8,0) and (-1.7,.5) .. (q1);
\draw[addflow] (1.7,0.8) .. controls (1.8,0) and (1.7,-.5) .. (q1);

\end{tikzpicture} 
    \caption{A Morse-Smale dynamical system on $S^2$ with a periodic orbit $O$}\label{fig2}
  \end{figure}
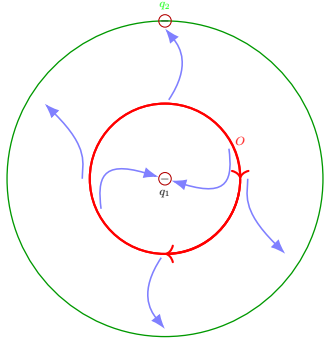

\begin{figure}
\centering
  \begin{tikzpicture}[scale=0.5, transform shape,
    every node/.style={font=\small},
    flow/.style={-Latex, thick},
    addflow/.style={-Latex, blue!50, line width=0.6pt},
    saddle/.style={circle, draw=black, inner sep=3.5pt},
    sink/.style={circle, draw=green!70!black, inner sep=0pt, minimum size=7pt},
    top/.style={circle, draw=black!70!black, inner sep=0pt, minimum size=7pt},
    closed/.style={red, thick},
  arrowcircle/.style={
    postaction={decorate},
    decoration={markings, mark=at position 0.25 with {\arrow{>}}}
  }
]

%--- outer disk (S^2 as disk model) ---
\draw[green!60!black, line width=.4pt] (0,0) circle (4.2cm);

%--- red closed orbit (γ) near top-left quadrant ---
\draw[closed] (0,0) circle (2cm);

\draw[<-,red, thick] (2,0) arc (0:350:2cm);

\draw[<-,red,thick] (-2,0) arc[start angle=180, end angle=0, radius=2cm];

%--- top (plus signs) ---
\node[top,label=above:$t_1$] (t1) at (0,2) {\(+\)};

%--- saddles ---
\node[saddle] (s1) at (0,-2) {};
\node[black] at (0.2,-1.8) {$s_1$};

%--- sinks (minus signs) ---
\node[sink,label=above:$q_1$] (q1) at (0,0) {\(-\)};
\node[sink,label=below:\green{$q_2$}] (q2) at (0,-4.2) {\(-\)};

%--- additional faint flow lines (blue) for visual context ---
\draw[addflow] (s1) .. controls (0,-3)  .. (0,-4);

\draw[addflow] (s1) .. controls (0,-1)  .. (q1);
\end{tikzpicture}
\caption{Replacing the periodic orbit in Fig. \ref{fig2} by two heteroclinic orbits, resulting in the situation of Fig. \ref{fig1}}\label{fig3}
\end{figure}
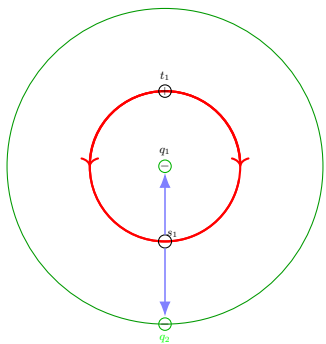
As in Fig. \ref{fig1}, in Fig. \ref{fig2} the outer green circle represents a single minimum point $q_2$. $q_1$ is likewise a minimum, but instead of a maximum and a saddle, we now have a critical closed orbit that is topologically equivalent to those two critical points. The orbit is unstable, and flow lines go towards the two minima. As shown in \cite{Eidi21, Eidi24}, we can still do a count of flow lines to construct a Floer type boundary that recovers the topology of the underlying space, $S^2$ in the present example. The count leads to the same homological structure as when we replace the periodic orbit by two heteroclinic orbits, in this case between a local maximum $t_1$ and a saddle $s_1$, as shown in Fig. \ref{fig3}. Again, $\partial_2 t_1=0, \partial_1 s_1=q_1+q_2, \partial_0=0$.
Analogously to Fig. \ref{fig1a}, in Fig. \ref{fig4}  (modified from \cite{Bannwart24}) counting flow lines does not yield the right boundary operator. Again, implicitly, there would be flow lines between critical objects of the same index. After turning the orbit to two fixed points and two heteroclinic orbits between them, the problem can be resolved as in Fig. \ref{fig1b}. Alternatively to preserve the orbit, we could first cancel one of the pairs $(s_1, q_1)$ or $(s_1, q_2)$, to obtain Fig. \ref{fig2}, in a process similar to the one allowed in the Morse cancellation theorem. The resulting system, as we call it "simple generalised Morse-Smale" in \cite{Eidi21}, would not have any flow lines between the elements of the same index and therefore a valid case to be handled via our generalised Floer-type boundary operator. Recall that, by the Morse cancellation theorem, we can cancel a pair of critical points $(p_k, p_{k-1})$ when there is a single flow line between them (if their unstable and stable manifolds, resp.,  intersect transversally). For more details on this theorem, we refer to \cite{Milnor, cancel}.

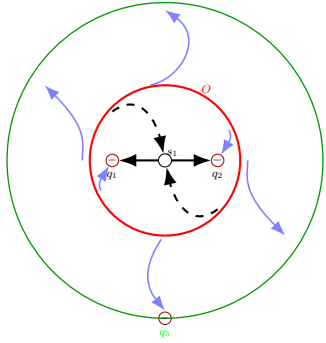
\begin{figure}
\centering
\begin{tikzpicture}[scale=0.5, transform shape,
    every node/.style={font=\small},
    flow/.style={-Latex, thick},
    addflow/.style={-Latex, blue!50, line width=0.6pt},
    saddle/.style={circle, draw=black, inner sep=3.5pt},
    sink/.style={circle, draw=red!70!black, inner sep=0pt, minimum size=7pt},
    closed/.style={red, thick}
]

%--- outer disk (S^2 as disk model) ---
\draw[green!60!black, line width=0.5pt] (0,0) circle (4.2cm);

%--- red closed orbit (γ) near top-left quadrant ---
\draw[closed] (0,0) circle (2cm);
\node[red] at (1.1,1.9) {$O$};

%--- saddles ---
\node[saddle] (s1) at (0,0) {};
\node[black] at (0.2,0.2) {$s_1$};

%--- sinks (minus signs) ---
\node[sink,label=below:$q_1$] (q1) at (-1.4,0) {\(-\)};
\node[sink,label=below:$q_2$] (q2) at (1.4,0) {\(-\)};

\node[sink,label=below:\green{$q_3$}] (q3) at (0,-4.2) {\(-\)};

%--- main separatrices (black) ---
\draw[flow] (s1) .. controls (-.7,0) .. (q1);
\draw[flow] (s1) .. controls (.7,0).. (q2);

%--- connections from closed orbit to saddles ---
\draw[flow,dashed] (-1.4,1.3) .. controls (-0.6,1.9) and (-0.2,0.9) .. (s1);
\draw[flow,dashed] (1.4,-1.3) .. controls (0.6,-1.9) and (0.2,-0.9) .. (s1);

%--- additional faint flow lines (blue) for visual context ---
\draw[addflow] (-2.2,0) .. controls (-2.2,0.4) and (-2.0,0.8) .. (-3.2,2);
\draw[addflow] (-.1,-2.1) .. controls (-.7,-3) and (-.4,-3.5) .. (0,-4);
\draw[addflow] (2.2,-0) .. controls (2.2,-.4) and (2,-.8) .. (3.2,-2);
\draw[addflow] (-1.7,-0.8) .. controls (-1.8,-.6) and (-1.7,-.5) .. (q1);
\draw[addflow] (1.7,0.8) .. controls (1.8,.6) and (1.7,.5) .. (q2);
\draw[addflow] (-.4,2) .. controls (.5,2.2) and (1,3.3) .. (0,4);
\end{tikzpicture}

    \caption{An example that does not satisfy the condition for the boundary operator}\label{fig4}
\end{figure}

\subsection{Combinatorial Morse theory}\label{morse}
Again, there is a discrete version, developed by Forman \cite{Forman98a,Forman98b,Forman98c}.
Here, a function assigns a value to every simplex or cell,
and certain inequalities between the values on a simplex and on its facets are required that can be
seen as analogues of the non-degeneracy conditions of Morse theory in the smooth setting.\\
We shall also use the presentation in \cite{Joharinad23}.
\begin{defi}\label{morsedef}
  A function $f:\Sigma \to \N$  is called a Morse function if for every
  simplex 
  $\sigma_k \in \Sigma_k$
  \begin{enumerate}
  \item there is at most one simplex $\rho_{k+1} \supset \sigma_k$
    with 
\bel{exch5} f(\rho_{k+1})\le f(\sigma_k)
\qe
and
\item at most one simplex $\tau_{k-1}\subset \sigma_k$ with
\bel{exch6} f(\tau_{k-1})\ge f(\sigma_k).
\qe
  \end{enumerate}
  The simplices for which neither of these possibilities holds are called \emph{critical}.
\end{defi}

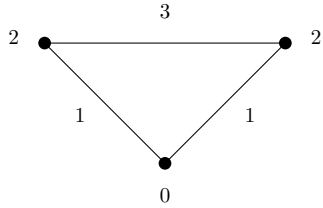
\begin{figure}
 \begin{center}
 \begin{tikzpicture}[scale =0.8, transform shape]
  \fill (6,3) circle (3pt);
  \fill (4,1) circle (3pt);
  \fill (2,3) circle (3pt);

   \draw (6,3) -- (4,1);
   \draw (2,3) -- (4,1);
   \draw (6,3) -- (2,3);

  \coordinate [label=above:3] (3) at ($ (4,3.3)$);
  \coordinate [label=left:2] (2) at ($ (1.7,3.1)$);
  \coordinate [label=right:2] (2) at ($ (6.3,3.1)$);
  \coordinate [label=left:1] (1) at ($ (2.8,1.8)$);
  \coordinate [label=right:1] (1) at ($ (5.2,1.8)$);
  \coordinate [label=below:0] (0) at ($ (4,.7)$);
  \end{tikzpicture}
 \end{center}

 \caption{A simplicial Morse function} \label{fig20}
 \end{figure}
  Fig. \ref{fig20} shows a Morse function with a critical vertex with value 0 (the minimum) and a critical edge with value 3 (the maximum). \\
    Again, we have 
\bel{morse11}
m_k\ge b_k\quad \text{ and } \sum_{k} (-1)^k b_k=\sum_{k} (-1)^k m_k
\qe
where $m_k$ now is the number of critical $k$-simplices.\\

For a Morse function, we draw an arrow from $\sigma_k$ to $\rho_{k+1}$ whenever $\sigma_k\subset
\rho_{k+1}$ and $f(\rho_{k+1})\le f(\sigma_k)$. For this arrow pattern, there is  at most one arrow per simplex,  because, as is easily shown,  for  a Morse function, it cannot happen that there  exist both a
simplex $\rho_{k+1}$ with 1. and a simplex $\tau_{k-1}$ with 2.\\ 
But not every such arrow pattern comes from  a Morse function. We need to exclude cycles, see \cite{Forman98b}. 
\\

Forman \cite{Forman98c} has also constructed a combinatorial version of Witten's approach described in Section \ref{witten}, conjugating the Eckmann Laplacian with $e^{-tf}$ for a Morse function. For the corresponding theory, he needs the additional condition on $f$ that when a $k$-simplex $\sigma$ is a facet of two different $(k+1)$-simplices $\tau^1,\tau^2 $, then $f(\sigma)\le \min (f(\tau^1),f(\tau^2))$. \\

From a Morse function, we can also construct a new, nonlocal Floer type \cite{Floer89} boundary  operator,  as first described in \cite{Forman98a}. For that purpose, we call a 
   simplex $\sigma_k$ \emph{upward noncritical} if there is a $\rho_{q+1}$ satisfying \rf{exch5}. We call it \emph{downward noncritical}  if there is a $\tau_{q-1}$ with \rf{exch6}. The  simplices satisfying neither of these  are the \emph{critical} ones. 

\begin{defi}
  Let $\sigma$  be a critical $k$-simplex of the discrete Morse function
$f$. Then for $\partial_F \sigma$, we select  $(k-1)$-simplices  by the following rules, with appropriate orientation rules. 
\begin{enumerate}
\item 
If $\rho \,< \,\sigma$ is critical then keep it.
\item If $\rho < \sigma $ is downward non-critical, ignore it.
\item 
If $\rho < \sigma$  is upward non-critical then there exists a unique
$\sigma ' \,>\, \rho$ \,with \,$f(\sigma ') \,\le \,f(\rho)$ (this
$\sigma'$ then is non-critical, hence in
particular different from $\sigma$). In this
case, apply these rules to $\sigma'$ in place of $\sigma$. That is, 
keep all the critical facets of this $\sigma '$, ignore the  downward
non-critical facets, and proceed with the  upward non-critical facets other than $\rho$
 as before. For any such $\rho'$, there has to be a unique $\sigma''\,>\, \rho'$ \,with \,$f(\sigma '') \,\le \,f(\rho')$, and since by 2., we have $f(\sigma ') \,> \,f(\rho')$, we conclude that  $\sigma'' \neq \sigma'$. Thus, we can go on. 
\end{enumerate}
\end{defi}
We have 
 \bel{fl2} \partial_F \circ \partial_F=0
\qe

The proof can proceed inductively, by inserting the arrows one by one. \\
A function that is constant on some  collections of simplices, but otherwise satisfies the rules of Def. \ref{morsedef}  is called a \emph{Morse-Bott function}, in analogy with the corresponding setting in the smooth case \cite{Yaptieu17}. As in the smooth case, we can also handle (finitely many) closed orbits of a simplicial flow \cite{Eidi21,Eidi24}. In fact, for a combinatorial vector field in Forman's terminology, a procedure  analogous to that of Franks \cite{Franks79} can be applied. Each closed orbit can be replaced by a pair of critical simplices whose dimensions differ by one, allowing us to construct a Floer-type boundary operator (as in the case of Morse-Smale dynamical systems) for these generalized dynamics in which the chain recurrent set consists not only of critical simplices but also of closed orbits.

\subsection{Discrete and continuous Morse theory}

{\bf Questions:} What are the relations between discrete and smooth Morse theory? Discrete Morse theory is essentially combinatorial while smooth Morse theory is geometric and analytic. Are there intermediates? For instance, if we consider a geometric realization of a  simplicial complex, can one translate a simplicial Morse function into a function on that geometric realization that behaves as a smooth Morse function? And as for the eigenvalues, what is the behavior under convergence of a simplicial complex or a graph to a Riemannian manifold?\\

In this section, we will provide some systematic answers. In fact, we shall describe three different approaches. 

So far, a simplicial complex $\Sigma$ has been treated as a combinatorial structure. We can also consider its geometric realization
\begin{equation}\label{geo1}
    |\Sigma|= \bigcup_{\sigma \in \Sigma} \conv\{\vec 1_{\{v\}}:v\in\sigma\}
\end{equation}
where {$\vec 1_{\{v\}}$ is the function that takes the value 1 on the vertex $v$ and 0 on all other vertices and $\conv$  denotes the convex hull. And so, $\conv\{\vec 1_{\{v\}}:v\in\sigma\}$ is the geometric simplex spanned by the vertices of $\sigma$. Of course, in $|\Sigma|$, all the face relations of $\Sigma$ are preserved.

There exists a theory of piecewise linear Morse functions on geometric simplicial complexes. A piecewise linear ({\bf PL} for short) function on a geometric simplicial complex is linear on each simplex and the linear pieces agree on common facets. Since a linear function has no critical points, a piecewise linear function on geometric simplicial complex can have critical points only at the vertices. So, this seems at odds with Forman's combinatorial concept where simplices of any dimension can be critical. The key to relating these two different concepts consists in going from a combinatorial \simp to a geometric one that is not the geometric realization of the former, but one whose vertices are the simplices of the combinatorial complex. This is the \emph{order complex} of $\Sigma$ whose simplices are the chains of simplices in $\Sigma$,
$$S_\Sigma:=\{\Cal \subset\Sigma: \Cal\text{ is a chain}\},$$
that is, a family of simplices ordered by inclusion. %\text{ any two members in }\C\text{ possess the inclusion relation}\}$$
%which collects all inclusion chains in $\K$.
Thus,  $\Cal$ is a {\sl chain} if for any $\sigma_1,\sigma_2\in\Cal$, either $\sigma_1\subset\sigma_2$ or $\sigma_2\subset\sigma_1$. 
The order complex  $S_\Sigma$  is itself a simplicial complex with vertex set
$\Sigma$.  According to \eqref{geo1}, the  geometric realization of $S_\Sigma$ is
$$|S_\Sigma|= \bigcup_{\Cal\in S_\Sigma}\conv\{\vec 1_{\sigma}:\sigma\in \Cal\}$$
where $\vec 1_{\sigma}$ is the function that takes the value 1 on the vertices of $\sigma$ and 0 on all other vertices. 
This geometric realization can be identified with the geometric realization of the barycentric subdivision $B_\Sigma$ of $\Sigma$. In that barycentric subdivision, also every simplex of $\Sigma$ becomes a vertex, as for $S_\Sigma$, and therefore, $S_\Sigma$ and $B_\Sigma$ are combinatorially the same, and hence their geometric realizations coincide. Also, this yields a digraph to which the method of path homology can be applied \cite{Grigoryan14}.

If $f$ then is a function on $\Sigma$ (to emphasize, a combinatorial function, that is, assigning a real number to every simplex of $\Sigma$), it becomes a function on the vertex set of $S_\Sigma$. And we then use the Lov\'asz extension from Section \ref{lovasz} to extend it as a function $f^L$ on the geometric realization 
$|S_\Sigma|$.

There is a general result \cite{Jost/Zhang24a}.
\begin{theorem}
\label{Mainthm:discrete-Morse}
For a simplicial complex with vertex set $V$ and face set $\Sigma$, let $f:\Sigma\to\R$ be an injective discrete Morse function. Then the following conditions are equivalent:
\begin{enumerate}
\item  $\sigma$ is a critical point of $f$   with $\dim \sigma=k$.
\item $\vec1_\sigma$ is a critical point of
  $f^L|_{|S_\Sigma|}$ with index $ k$ in the sense of metric Morse theory.
  \item  $\vec1_\sigma$ is a Morse critical point of $f^L|_{|S_\Sigma|}$ with index $ k$ in the sense of topological Morse theory.
\item  $\vec1_\sigma$ is a critical point of $f^L|_{|S_\Sigma|}$ with index $ k$ in the sense of PL Morse theory. 

\end{enumerate}
Moreover, the discrete Morse vector $(n_0,n_1,\cdots,n_d)$,  representing the number $n_i$ of critical points with index $i$, of $f$ coincides with the  continuous Morse vector of $f^L|_{|S_\Sigma|}$.

Conversely, if the geometric simplicial complex $|S_\Sigma|$ admits a Morse vector $\vec c$, then there exists a simplicial complex $\Sigma'$ that is PL-homeomorphic to $\Sigma$ and has the (discrete) Morse vector $\vec c$.

\end{theorem}

In order to explain this theorem, we need to describe the concepts involved. But let us first note that we don't assume the geometric complex $|\Sigma|$ is a topological manifold, and therefore  topological results on manifolds like \cite{Gallais10,Benedetti10,Benedetti16} do not apply. (For instance,  Benedetti \cite{Benedetti10} showed that  if a PL triangulation of a manifold admits a smooth Morse function with 
$m_k$ critical points of index $k$, then after an appropriate subdivision, one can construct a boundary-critical discrete Morse function whose interior critical faces correspond exactly to these smooth critical points, with 
$m_k$ faces of dimension $k$.)

\begin{enumerate}
\item \textbf{Discrete Morse theory} has already been described above in Section \ref{morse}. 
\item \textbf{Metric Morse theory}:
Let $F:M\to \R$ be   a continuous function on the metric space $(M,d)$. For $p\in M$, 
$$U(p,\delta)=\{q\in M: d(p,q)<\delta\}$$
is the open ball in $M$ of radius $\delta$ about $p$. The {\sl weak slope} \cite{Degiovanni94,Ioffe96,Katriel94}  $|dF|(p)$ is defined as the supremum of  $\epsilon\ge 0$ for which 
 there exist $\delta>0$ and a continuous map
$$H:U(p,\delta)\times[0,\delta]\to M$$
satisfying $$F(H(q,t))\le F(q)-\epsilon t,\;\; d(H(q,t),q)\le t$$
for all $q\in U(p,\delta)$ and $t\in[0,\delta]$. Points $p$ with vanishing weak slope, i.e., $|dF|(q)=0$ are called \emph{critical} for $F$. These are the points for which no neighborhood can be deformed entirely to smaller values of $F$. 

The local behaviour of $F$ near a critical $p$ is captured  by the {\sl
  critical group} 
  $$C_*(F,p):=H_*(\{F\le c\}\cap U(p,\delta),\{F\le c\}\cap
U(p,\delta)\setminus \{p\};\R),$$
for suitable $\delta >0$. Here $H_*(\cdot,\cdot;\R)$
 is the singular relative homology. In general, $C_k(F,p)$ may be non-trivial for several $k\in \Z$, as we do not impose a Morse type assumption on $F$ here. 
 We then obtain the Morse polynomial $\sum_{k}\mathrm{rank}\, C_k(F,p) t^k$. And the indices (possibly several) of a critical point $p$ are the values of $k$ for which the coefficients are non-zero. 

\item \textbf{Topological Morse theory}:
Let $F:M\to \R$ be   a continuous function on the 
 topological space $M$. $p\in M$ is a {\sl  regular point} of $F$ if there exists a neighborhood $U(p)$  and a continuous map $$H:U(p)\times[0,1]\to M,\;\; H(q,0)=q$$
with  $$F(H(q,t))< F(q),$$
for all $q\in U$ and $t>0$. Points that are not regular are called  {\sl critical} for  $F$ on $M$. Critical groups and indices can be defined as in the metric setting \cite{Degiovanni10}.

\item \textbf{Piecewise-Linear Morse theory}: This theory works on a geometric simplicial complex $M$ considered as a PL manifold.  We will adopt the definitions of  K\"uhnel \cite{Brehm87,Kuhnel90} and  Edelsbrunner  and Harer  \cite{Edelsbrunner10}. \\
Let $F:M\to \R$ be   a PL  function. As already mentioned, we need to pay attention only to the behavior of $F$ near the vertices $v$ of $M$. 
     $\mathrm{star}_-(v)$ then indicates the subset of the star of $v$ on which  $F(q)\le F(v)$, and analogously $\mathrm{link}_-(v)$. 
     Here, $\mathrm{star}(v)$  is the subcomplex consisting of all simplices of $M$ that contain $v$, together with all their faces while $\mathrm{link}(v)$ is the subcomplex consisting of all simplices that are disjoint from $v$ but co-appear with $v$ in some higher-dimensional simplex of $M$. The latter are those simplices in $\mathrm{star}(v)$ that do not contain $v$.
\begin{defi}
A vertex $v$ of $M$ is  a PL critical point of $F$ with index $k$ and multiplicity $m_k$ if the Betti number $b_k(\overline{\mathrm{star}_-(v)},\mathrm{link}_-(v))=m_k$.
\end{defi}
Alternatively, and equivalently, we can use the rank of the reduced $(k-1)$-th homology group of $\mathrm{link}_-(v)$. 
Clearly, a PL critical point may have many indices and multiplicities. The PL  function $F$ is called a PL Morse function if all critical vertices $v$ are non-degenerate in the sense that  $\sum_k m_k(v)=1$.
\end{enumerate}

\begin{remark}
As a bridge between discrete Morse theory and discrete Laplacian on simplicial (cell) complexes, Forman in \cite{Forman98c} presented a discrete counterpart to Witten’s Hodge-theoretic approach \cite{Witten82} to Morse theory, as already mentioned above. Inspired by Witten’s analytic proof of the classical Morse inequalities via deformed Laplacians, Forman constructed a combinatorial version that mirrors the key structures of the smooth theory while avoiding its analytical difficulties, such as infinite-dimensional analysis and transversality issues. His framework establishes combinatorial analogues of Witten’s main ingredients, including the formulation of a discrete Hodge theory and a combinatorial proof of the Morse inequalities. In fact, Forman’s theory applies  to arbitrary cell complexes, not only to those arising from cell decompositions of manifolds, thereby extending Witten’s ideas from global analysis to a purely combinatorial realm. \end{remark}

\begin{remark}
The authors \cite{Jost/Zhang24a} have also developed a discrete Morse theory for certain types of hypergraphs, in which a topological lemma concerning general hypergraphs \cite[Proposition 5.6]{Jost/Zhang24a} can be regarded as a hypergraphic version of a result of Lickorish \cite[Lemma 1]{Lickorish}.
\end{remark}

\section{Cheeger type inequalities}
\subsection{Cheeger cuts}\label{cheeger}
\subsubsection{For Riemannian manifolds}
As before, we consider a compact and connected Riemannian manifold $M$. 
We order the eigenvalues of the Laplace-Beltrami operator of  $M$ as
\begin{equation}
    0=\lambda_1 <\lambda_2\le \lambda_3 \le\cdots
\end{equation}
where the eigenfunction for $\lambda_1 =0$ is constant and $\lambda_2>0$ since $M$ is connected, that is, $b_0=1$. 
Jeff Cheeger \cite{Cheeger70} used $\lambda_2$ to control a
  global quantity,
\bel{cheeger1}
h(M):= \inf \frac{\Vol_{d-1}(S)}{\min(\Vol_d(M_1),\Vol_d(M_2))}
\ee
where we consider  all $(d-1)$-dimensional submanifolds, i.e., hypersurfaces, 
$S$  that divide $M$ into two pieces $M_1,M_2$, with their volume $\Vol_{d-1}$. $h(M)$ is small when $M$ can be divided into two  pieces $M_1$
and $M_2$ of large volume by
cutting $M$ along a hypersurface $S$ of small volume. Here, $\Vol_{d-1}(S)$ has to be taken as the perimeter of $S$, as we do not know a priori whether $S$ is smooth, and for the infimum in \eqref{cheeger1}, we need to admit all hypersurfaces  of finite perimeter. In other words, we have to consider the BV functional.  That is,
\begin{equation}\label{bv1}
  \Vol_{d-1}(S)=\|D\chi_{M_1}\|_{BV} =\|D\chi_{M_2}\|_{BV}
\end{equation}
where $\chi_U$  denotes the characteristic function of the set $U$. For a brief introduction to the BV functional, see for instance \cite{Jost98}. (There also is the issue of the regularity of a minimizing hypersurface $S$, but this is not addressed here.)\\

Cheeger's estimate says
\bel{cheeger2}
\lambda_2 \ge \frac{h(M)^2}{4}.
\ee
Buser \cite{Buser78,Buser80} showed that this estimate cannot be improved.\\

Let us sketch the argument for \eqref{cheeger2}.  An
eigenfunction $u_2$ for the eigenvalue $\lambda_2$ satisfies
$\int_Mu_2=0$, as 
$u_2$ must be orthogonal to the $\lambda_1$-eigenspace (the constants), and therefore divides  $M$ into $M_1,M_2$ with $u_{2}|_{M_1}<0, u_{2}|_{M_2}>0$;  we may assume 
\begin{equation}\label{cheeger2a}
    \Vol_d(M_1)\le \Vol_d(M_2).
    \end{equation}
    We multiply the eigenvalue equation 
\bel{cheeger3}
\Delta u_2=\lambda_2 u_2 \text{ in } M_1,\qquad  u_2=0 \text{ on }
\partial M_1
\ee
 by $u_2$ and integrate by parts over $M_1$, suppressing the technical issue of 
the regularity of  $\partial M_1$. This yields 
\bel{cheeger4}
\lambda_2=\frac{\int_{M_1} \langle du_2, du_2 \rangle}{\int_{M_1}
  \langle u_2, u_2 \rangle}.
\ee
 $\chi:=(u_2)^2 (=\langle u_2,u_2\rangle)$ satisfies $d\chi =2u_2
du_2$, and so 
\bel{cheeger5}
\left(\int_{M_1} \| d\chi \|\right)^2 \le 4 \int_{M_1} \langle du_2, du_2
\rangle\ \int_{M_1}
  \langle u_2, u_2 \rangle. 
\ee
From \eqref{cheeger4} then 
\bel{cheeger6}
\lambda_2 \ge \frac{1}{4}\left(\frac{\int_{M_1} \| d\chi
  \|}{\int_{M_1}\chi}\right)^2.
\ee
With the  coarea formula 
\ba
\nonumber
\int_{M_1} \|d\chi\| &=&\int_{-\infty}^{\infty}
\Vol_{d-1}(\chi^{-1}(t)\cap M_1) dt\\
\label{cheeger8}
&\ge& h(M) \int_0^\infty \Vol_d(\{\chi\ge
t\}\cap M_1)dt \\
\label{cheeger8a}
&=&h(M)\int_{M_1}\chi.
\ea
 \rf{cheeger8} is justified since for $t>0$, $\chi^{-1}(t)$ is a closed
hypersurface that does not meet $\partial
M_1$. And $\Vol_d (\{\chi\ge
t\})\le \Vol_d (M_1)\le \Vol_d (M_2)$, by \eqref{cheeger2a}.\\
\rf{cheeger2} then follows from \rf{cheeger8a} and \rf{cheeger6}.\\

Buser \cite{Buser82} showed that $h(M)$ can also be controlled 
from below by $\lambda_2$,
\begin{equation}\label{cheeger10}
    \lambda_2 \le 2a (n-1)h(M) +10 h(M)^2
\end{equation}
if the Ricci curvature of $M$ is bounded from below by $-(n-1)a^2$. In particular, when $a=0$, that is, when $M$ has non-negative Ricci curvature, Buser's inequality says
\begin{equation}\label{cheeger11}
    \lambda_2 \le 10 h(M)^2
\end{equation}

\subsubsection{For graphs}\label{graph}
A Cheeger type inequality for graphs was first obtained by Dodziuk \cite{Dodziuk84} and Alon-Milman \cite{Alon85}. The version of the inequality \eqref{che13a},\eqref{che13b} that we shall use below is due to Chung \cite{Chung97}. Our exposition also partly follows \cite{Jost/Mulas/Zhang26}. 

Let $\Gamma$ be a graph with vertex set $V$ and edge set $S$. The edges $e=(v_1,v_2)$ may have symmetric  non-negative weights
\begin{equation*}
    w_{v_1,v_2}=w_{v_2,v_1}\ge 0
\end{equation*}
where the value 0 simply indicates absence of the corresponding edge. Unless stated otherwise, we assume that the graph is connected, that is, any two vertices can be connected by a chain of edges of positive weights. And unless stated otherwise, we assume that the graph has no loops, that is, for all edges $e=(v_1,v_2)$, $v_1\neq v_2$. \\
We write $u\sim v$ if $w_{u,v}>0$. The degree of a vertex is
\begin{equation*}
    \deg v=\sum_{u\sim v}w_{u,v},
\end{equation*}
and the volume of 
a subset $S$ of the vertex set $V$ is
\begin{equation*}
    \Vol S =\sum_{v\in S}\deg v
\end{equation*} and its complement is
$\overline{S}=V\setminus S$. 
For  $V_1,V_2\subset V$,
\begin{equation*}
 |E(V_1,V_2)|=\sum_{u\in V_1, v\in V_2}w_{u,v} .  
\end{equation*}
When the graph is unweighted, all edges that are present have the weight 1, and $\deg v$ counts the number of vertices connected with $v$ by an edge, and $|E(V_1,V_2)|$ is the number of edges between $V_1$ and $V_2$. 

We then put
\bel{che17a}
\eta(S)= \frac{|E(S,\overline{S})|}{\min(\mathrm{vol}(S),\mathrm{vol}(\overline{S}))},
\qe
and 
\bel{che17b}
h=  \min_{S}\eta(S).
\qe
This constant was first considered by Polya \cite{Polya51}, but it is now called the Cheeger constant, because it is analogous to the Cheeger constant \eqref{cheeger1} in Riemannian geometry, and an estimate analogous to \eqref{cheeger2} obtains, 
\begin{equation}\label{che13b}
\frac{1}{2} h^2\le \lambda_2. 
\end{equation}
In fact, one has the slightly stronger estimate, and also an easier upper bound,  
\begin{equation}\label{che13a} 1-\sqrt{1- h^2}\le  \lambda_2 \le 2 h.
\end{equation}
Here, $\lambda_2$ is the second eigenvalue, that is, the first non-vanishing one of the normalized graph Laplacian which we shall now introduce. This is a special case of the Eckmann Laplacian $L_0=L_{0}^{up}=\delta_0^* \circ \delta_0$ \eqref{eck14}, \eqref{eck12}, but to make it concrete we need to specify the scalar product used for the definition of $\delta_0^*$. For two functions $f,g$ on the vertex set $V$, it is given by
\begin{equation}\label{ch41}
    (f,g)=\sum_v \deg v f(v) g(v).
\end{equation}
We then have
\begin{equation}\label{che42}
 \Delta f(v):=L_0    f(v) = f(v)-\frac{1}{\deg v}\sum_{u\sim v}w_{uv} f(u)
\end{equation}
and this is the normalized graph Laplacian for which \eqref{che13a} holds.\\

Without the factor $\deg v$, we obtain the algebraic graph Laplacian 
\begin{equation}\label{che43}
 A f(v)=\deg v f(v)-\sum_{u\sim v}w_{uv} f(u)
\end{equation}
for which also Cheeger type inequalities hold, see \cite{Dodziuk84,Alon85}.

\begin{remark}
There seems to be some inconsistency here. In \eqref{che13b}, the lower bound is $\frac{1}{2}h^2$, whereas in \eqref{cheeger2}, we have $\frac{1}{4}h^2$. And indeed, there are inconsistent conventions that explain the contrast between the factors $\frac{1}{2}$ and $\frac{1}{4}$. For graphs, we have  counted each edge only once in the total variation functional and in the numerator of the Rayleigh quotient.  But in  the setting of Riemannian geometry, we  count all vertices and the edges connected to each vertex, and consequently, we count each edge twice. If we used such a counting rule, we should redefine the Cheeger constant of a graph as $2h$, and consider $2\lambda_2$ as  the first nontrivial eigenvalue of the graph Laplacian. Then the inequality \eqref{che13b} becomes $2\lambda_2\ge\frac{1}{4}(2h)^2$ which agrees with the manifold Cheeger inequality \eqref{cheeger2}.
\end{remark}
\begin{remark}
As stated in Section \ref{plap}, the $p$-Laplacian Cheeger inequality refines \eqref{che13b}. By the monotonicity property that $p(2\lambda_k(\Delta_p))^{\frac1p}$ is increasing with respect to $p\in[1,+\infty)$ (see \cite{Zhang25}) and the Cheeger identity $h=\lambda_2(\Delta_1)$ (see \cite{Hein10}),   
 %\frac{2^{p-1}}{p^p}h^p\le \lambda_k(\Delta_p)
$\frac{p}{2}(2\lambda_2(\Delta_p))^{\frac1p}$ decreasingly converges to $h$ as $p$ approaches $1$ from the right, and in particular, when $p=2$, we have $\sqrt{2\lambda_2(\Delta_2)}\ge h$ which reduces to the Cheeger inequality \eqref{che13b}.
\end{remark}

\subsubsection{Signed graphs and dual Cheeger inequalities}\label{dual}
In this section,  $\Gamma$ is a graph with $N$ vertices that is  unweighted, or equivalently, any edge $e=(u,v)$ that is present has weight 1.   But here, the graph in addition carries a sign
function on the edges. That is, we consider a  \emph{signed graph}, a graph $\Gamma$  with a map $s$ from its edge set
to $\pm 1$.  Signed graphs were first introduced and studied by Harary \cite{Harary55,Harary57,Harary59,Harary60}. The systematic mathematical study and the derivation of Cheeger inequalities was achieved by Atay and Liu \cite{Atay20}, which is our main reference for this section.\\
When  we change the signs of all edges containing a particular vertex, we speak of a \emph{sign switch}. We shall see that such sign switches do not affect the properties that we are interested in. Therefore, we can try to bring the signs into a particular configuration by switching the signs of appropriate vertices. The best one can hope for is that the graph thereby becomes \emph{balanced}, that is, 
all signs $=1$, or  \emph{antibalanced}, all signs $=-1$. (These two cases are not mutually exclusive.) In particular, we can equip every edge of an ordinary unweighted graph without signs with the sign $+1$ and therefore consider it as a balanced graph.

%Using the framework of signed graphs developed  in Section \ref{signed},  
We define the   Laplacian  of the signed graph $(\Gamma,s)$ as
\begin{equation}
  \label{slap}
  \Delta_s f(v)= f(v)-\frac{1}{\deg v}\sum_{v' \sim v}s(vv')f(v')=\frac{1}{\deg v}\sum_{v' \sim v}(f(v)-s(vv')f(v')).
\end{equation}
One checks that 
\begin{enumerate}
    \item The spectrum $\lambda_1\le \lambda_2 \le \cdots\le \lambda_N$ of $\Delta_s$ is contained  in the interval $[0,2]$. 
    \item  $\lambda_1=0$ iff $(\Gamma,s)$ is balanced.
    \item $\lambda_N=2$ iff $(\Gamma,s)$ is antibalanced. 
\end{enumerate}
As a consequence of 2. and 3., an ordinary graph without signs, which, as explained above, can be seen as a balanced graph, is also antibalanced iff bipartite.\\

For  disjoint $V_1,V_2\subset V$, $$E^+(V_1,V_2)=\{\{u,v\}\in E:u\in V_1,v\in V_2,s(uv)=1\}
$$ and
$$
E^-(V_1)=\{\{u,v\}\in E:u,v\in V_1,s(uv)=-1\}$$ 
and define the \emph{signed bipartiteness ratio} as
\begin{equation}\label{bip}
\beta^s(V_1,V_2)=\frac{2
\left(|E^-(V_1)|+|E^-(V_2)|+|E^+(V_1,V_2)|\right)+|\partial(V_1\sqcup V_2)|}{\vol(V_1\sqcup V_2)}.
\end{equation}
The  \emph{signed Cheeger constant} of  $(\Gamma,s)$ is  
$$h_1^s=\min\limits_{(V_1,V_2)\ne(\emptyset,\emptyset)}\beta^s(V_1,V_2)$$
where $V_1$ and $V_2$ should be disjoint, but we do not require $V_1\cup V_2=V$.\\
Importantly, $\beta^s$ and  $h^s$ are switching invariant.  

We have the  Cheeger inequality for signed graphs of  \cite{Atay20},
\begin{equation}
  \label{atay}
1-\sqrt{1-(h_1^s)^2} \le \lambda_1(\Delta_s)\le 2h_1^s.  
\end{equation}
Let, as before $V_1\cap V_2=\emptyset $ and put $\overline{V_1\cup V_2} =: V_3$ to obtain  a
partition of the vertex set $V$. 
We then have the
\emph{dual Cheeger constant}\index{dual Cheeger constant} of \cite{Bauer13}
\be\label{hh}\overline{h}(\Gamma)=  \max_{V_1,V_2}
\frac{2|E(V_1,V_2)|}{\mbox{vol}(V_1)+\mbox{vol}(V_2)}.\qe
Alternatively, we have the 
  \emph{bipartiteness ratio} of \cite{Trevisan12}
  \bel{dual2}
  \beta(\Gamma)=\min_{V_1,V_2} \frac{2|E(V_1,V_1)| +2|E(V_2,V_2)| +|E(V_1
    \cup V_2, \overline{V_1 \cup V_2})|}{\mbox{vol}(V_1)+\mbox{vol}(V_2)}.
  \qe
 Since $\vol S= 2 |E(S,S)|+ |E(S,\overline{S})|$, 
  \begin{equation}\label{dual3}
    \beta(\Gamma)=1-\overline{h}(\Gamma),
  \end{equation}
  and since we have
  \begin{equation}
      \beta(\Gamma)=h_1^{\sigma_-}
  \end{equation}
where $\sigma_-\equiv -1$, we obtain a  \emph{dual Cheeger inequality} controlling  $\overline{h}$ and the largest eigenvalue $\lambda_N$ of an (unsigned) graph in terms of each other, 
  \begin{equation}\label{A2} 2\overline{h} \leq
\lambda_{N} \leq 1+\sqrt{1-(1-\overline{h})^2}.\end{equation}
This is a result of \cite{Bauer13}. The analogy with \eqref{che13a} is obvious, and it is best understood in terms of signed graphs, even though the graph to which it is applied was unsigned. \\
 
The $j$-way signed Cheeger constant of \cite{Lee12,Atay20} is defined as
$$h_j^s=\min\limits_{\{(V_{2i-1},V_{2i})\}_{i=1}^j}\max\limits_{1\le i\le j}\beta^s(V_{2i-1},V_{2i})$$
where we allow for  all possible $j$  pairs of
disjoint sub-bipartitions $(V_1,V_2)$, $(V_3,V_4)$, $\ldots$,
$(V_{2j-1},V_{2j})$.  $h^s_j$ is again switching invariant.  

We then have the higher order Cheeger inequalities of Atay and Liu \cite{Atay20} which generalize the  inequalities of \cite{Lee12},
\begin{equation}\label{bip3}\frac{(h_j^s)^2}{Cj^6}\le \lambda_j(\Delta_s)\le 2h_j^s
\end{equation}
 for all
signed graphs $(\Gamma,s)$, and all $j\in\{1,\ldots,N\}$, with
an absolute constant $C$.

\subsubsection{From graphs to manifolds}
When we construct a graph by random sampling from a compact Riemannian manifold $M$ and connecting points at distance $\le \epsilon$, with $\epsilon \to 0$ at an appropriate rate as the number $n$ of sample points increases, then the suitably rescaled spectra of the normalized graph Laplacians converge to the spectrum of the Laplace-Beltrami operator on $M$, as first exploited in \cite{Belkin03}, with quantitative estimates provided in \cite{Garcia20}. \\
Also, it was shown in \cite{Garcia16} that the perimeter functionals of the graphs $\Gamma$-converge to that of $M$. $\Gamma$-convergence is a notion of variational convergence for functionals, introduced by De Giorgi \cite{DeGiorgi75}, implying that the minimizers of approximating functionals converge to those of the limit functional. For a short introduction, see for instance \cite{Jost98}. In \cite{Garcia16}, it is shown that the BV functionals on the approximating graphs when functions on the graph are extended to piecewise constant functions on the manifold, with local  values coming from those on the graph vertices, converge to the BV functional on the Riemannian manifold.

\subsection{Higher dimensional Cheeger relations}\label{higher}
\subsubsection{Signed graphs constructed from simplicial complexes}\label{sign}
Since we work with orientations of simplices, that is, with signs, it is natural that signed graphs can play a role for encoding some of the combinatorial structure of a simplicial complex. \\
As we have explained in Section \ref{simp}, the Laplacians depend on the choice of scalar products on the cochains. And in Section \ref{graph}, we have used a particular such product to define the normalized graph Laplacian. We extend that definition here to higher dimensional cochains, to obtain normalized Laplacians for which Cheeger inequalities will hold.

The  normalized up-Laplacian of a  simplicial complex $\Sigma$ is
\begin{equation}
\label{asc1}
(\Delta_{k}^{up} f)([\sigma])=  f([\sigma]) 
-  \frac{1}{\deg \sigma}\sum_{\substack{\sigma'\in \Sigma_{k}: \sigma'\neq \sigma,\\ \exists\tau\in\Sigma_{k+1}:\sigma,\sigma'\in \partial \tau}} s([\sigma],[\sigma'])f([\sigma']),    
\end{equation}
where 
\begin{align}
s([\sigma],[\sigma'])&:= 
-\sgn([\sigma],\partial [\tau])\sgn([\sigma'],\partial [\tau])\label{asc2}
%\ea        
\end{align}
and  $\sgn=+1 (-1)$ when the two orientations agree (disagree). 

Thus, $s$ yields a sign, and  we can consider the  signed graph $(\Gamma_k,s)$ whose vertices are  the $k$-simplices of our simplicial complex, and we connect two
different such vertices $\sigma, \sigma'$  by an edge, $\sigma
\sim \sigma'$, if they are both facets of the same $(k+1)$-simplex. \\
$\Delta_{d}^{up}$ can be expressed through  the Laplacian
$\Delta_{(\Gamma_k,s)}$ of this signed graph. But the two are not identical, since the degree of $k$-simplex in our complex is different from its degree when we consider it as a vertex in $\Gamma_k$. In fact, when a $k$-simplex is the facet of a $(k+1)$-simplex, that contributes 1 to its degree, but in $\Gamma_k$, it is then connected with the $k+1$ other facets of that simplex, and so here we get a contribution of $k+1$ to the vertex degree. 
Thus 
\begin{equation}
 \label{asc3}
\Delta_{k}^{up}= (k+1)\Delta_{(\Gamma_k,s)} - k\ \mathrm{Id}.
\end{equation} 
Therefore, the eigenvalues $\mu_j$ of $\Delta_{k}^{up}$ can be computed from  the eigenvalues $\lambda_j$ of $\Delta_{(\Gamma_k,s)}$  as
\begin{equation}
    \label{asc5}
\mu_j =(k+1)\lambda_j -k.
\end{equation}
Since, as observed in Section \ref{dual}, the eigenvalues of $\Delta_{(\Gamma_k,s)}$ satisfy $0\le \lambda_j\le 2 $,
those of 
   $\Delta_{k}^{up}$ satisfy $0\le \lambda_j\le k+2 $.
   But this also constrains the eigenvalues of $\Delta_{(\Gamma_k,s)}$. 
Since $\mu_j \ge 0$, they are $\ge  \frac{d}{d+1}$, with equality
 iff there is some non-trivial $f$ with $\delta_k
  f=0$, and  the multiplicity of this eigenvalue $\frac{k}{k+1}$ then is  the dimension of the kernel of the coboundary
  operator $\delta_k$, that is, the Betti number $b_k$ of $\Sigma$.  
In particular, $(\Gamma_k,s)$ is never balanced for $k>0$, according to the criterion observed in Section \ref{dual}.\\
From the considerations in Section \ref{dual}  and the relation \eqref{asc5}, we also conclude that 
the spectrum of $\Delta_{k}^{up}$ contains the  eigenvalue $k+2$ iff $(\Gamma_k,s)$ possesses an antibalanced component, and the multiplicity of  this eigenvalue $k+2$ is  number of antibalanced components of $(\Gamma_k,s)$.

We consider the   signed graph $(\Gamma_k,-s)$ with the opposite sign as $(\Gamma_k,s)$, with its Laplacian  $\Delta_{(\Gamma_k,-s)}$ with the eigenvalues $\lambda_1^-\le \dots \le \lambda_M^-$, then instead of \eqref{asc5}, we have
\begin{equation}
    \label{asc6}
\mu_j =k+2 -(k+1)\lambda_{M-j}^-.
\end{equation}
Summarizing \eqref{asc5} and \eqref{asc6}, $\mu$ is an eigenvalue of $\Delta_{k}^{up}$ iff
$\frac{\mu +k}{k+1} $ is an eigenvalue of $\Delta_{(\Gamma_k,s)}$ iff $\frac{k+2-\mu}{k+1} $ is an eigenvalue of
$\Delta_{(\Gamma_k,-s)}$.
And the multiplicity of the eigenvalue $0$ of $\Delta_{k}^{up}$ is
$\ge k+1$,  
while the multiplicity of the eigenvalue $k+2$ of $\Delta_{k}^{up}$ is equal to 
the number of balanced components of $\Delta_{(\Gamma_k,-s)}$.\\
We can now turn to relating the smallest positive and the largest eigenvalue of $\Delta_{k}^{up}$ to geometric quantities. 

\subsubsection{Disorientability}\label{disorient}
From \cite{Rosenthal17}, we recall that a \simp is disorientable if we can orient the simplices of the maximal dimension $n$ in such a way that they always induce the same orientation on a joint facet. \\

A graph is disorientable iff it is bipartite iff it has no cycles of odd length iff its largest eigenvalue is $=2$, see Lemma \ref{dis} below. For generalizing this to higher dimensional simplicial complexes, we consider the signed graph $(\Gamma_{n-1},-s)$ of Section \ref{sign} . Its vertices are the $(n-1)$-simplices, and when two such simplices are facets of the same $n$-simplex $\tau$, then $-s=+1$ when the orientation induced by  $\tau$ is either compatible or incompatible for both of them. And $-s=-1$ if precisely one of the induced orientations is compatible.  Therefore, for $n\ge 2$, the graph $(\Gamma_{n-1},-s)$ can be balanced only if all induced orientations are compatible, because we cannot have three or more incompatible signs. Therefore, if the graph is balanced, every $n$-simplex has to introduce compatible orientations on all its facets. And therefore, conversely, two $n$-simplices need to introduce the same orientation on joint facets. That is, the \simp has to be disorientable. And $(\Gamma_{n-1},-s)$ is balanced precisely if it has the eigenvalue 0. And by \eqref{asc6}, this is equivalent to $\Delta_{n-1}^{up}$ having the eigenvalue $n+1$. So, we conclude
\begin{lemma}\label{dis}
  A \simp is disorientable iff its normalized up-Laplacian $\Delta_{n-1}^{up}$ achieves the maximal possible eigenvalue $n+1$.   
\end{lemma}
\qed

As this maximal eigenvalue is $=2$ for $n=1$, that is, a graph, this generalizes the result for  graphs just mentioned. 
\begin{remark}
A proof of  Lemma \ref{dis} based on direct Laplacian computations is provided in  \cite[Proposition 2.7]{Rosenthal17}. 
\end{remark}

The question of a geometric characterization 
 in higher dimensions was clarified in \cite{Eidi24a,Eidi25}, as we now want to describe. We need some definitions. A simplex is called \emph{branching} if its degree is $>2$ (and so, orientable simplicial complexes cannot have branching $(n-1)$-simplices). A $k$-cycle of length $\ell$ consists of $k$-simplices $\sigma_0,\dots ,\sigma_\ell$
with $\sigma_\ell=\sigma_0$ so that $\sigma_j$ and $\sigma_{j+1}$ always have a common cofacet, that is, for each $j$ there exists a $(k+1)$-simplex $\tau_j$ of which both $\sigma_j$ and $\sigma_{j+1}$ are facets. A cycle is a subcomplex of $\Sigma$, and therefore a \simp in its own right. A non-orientable cycle is called \emph{twisted}. For $0\le k\le n-1$, we can construct an up-dual graph 
$\Gamma_k^{up}$ of $\Sigma_k$, and for $1\le k\le n$ a down-dual graph $\Gamma_k^{down}$. In either case, the vertices are the $k$-simplices, and in $\Gamma_k^{up}$, $\sigma_k \sim \sigma'_k$ (i.e., they are connected by an edge) iff they have a common cofacet $\tau$ (which then is unique), while in $\Gamma_k^{down}$, $\sigma_k \sim \sigma'_k$ if they share a facet (which likewise will be unique). Later on, we will consider $\Gamma_k^{up}$ for oriented simplices $[\sigma_k], [\sigma'_k]$ and equip each edge with a sign, depending on whether their cofacet $[\tau_{k+1}]$ induces different or identical  orientations on them. But for the moment, we look down, at $\Gamma_n^{down}$ where $n=\dim \Sigma$. It is shown in \cite{Eidi24a} that $\Sigma$ is disorientable if the only odd cycles in $\Gamma_n^{down}$ come from branching $(n-1)$-simplices or twisted cycles in $\Sigma_{n-1}$, and there are no twisted cycles of even length. (For the normalized Laplacian, this was shown in \cite{Horak13} by a different argument.)\\
In Section \ref{largest} below, we shall obtain another characterization for the largest eigenvalue of $\Delta_{n-1}^{up}$ to be $n+1$, namely that the Cheeger constant $h_1(\Sigma_{n-1})$ vanishes, see \eqref{bip4}.

\subsubsection{The smallest positive eigenvalue}\label{smallest}
Above, we have discussed the Hodge and Eckmann Laplacians as generalizations of the Laplace-Beltrami operator and the graph Laplacian. So, naturally the question arises what geometric information their eigenvalues encode. The eigenvalue 0 contains topological information. Its multiplicity yields the Betti number of the corresponding dimension. But what about the non-vanishing eigenvalues? In particular, does the smallest non-vanishing eigenvalue tell us anything about higher order decompositions of the underlying object, analogous to the Cheeger inequality? 

This problem was initiated in \cite{Dotterrer12,Steenbergen14}. A solution was given in \cite{Jost/Zhang24b}. In particular, this solution provides an affirmative answer to the first part of a question posed by Dotterrer and Kahle \cite[Question 1]{Dotterrer12}. We shall now describe that result. It applies to the up-Laplacian, but since the non-zero spectrum of the down-Laplacian $\Delta_k^{down}$ equals that of the up-Laplacian $\Delta_{k-1}^{up}$, it also applies to the former, although the geometric interpretation will be in terms of the up-Laplacian.\\

 We recall the \emph{normalized} Laplacian
 $\Delta_{k}^{up}$, which according to \eqref{asc1} is given by
 \begin{equation}\label{nor1}(\Delta_{k}^{up} f)([\sigma])=  f([\sigma]) 
+ \frac{1}{\deg \sigma}\sum_{\substack{\sigma'\in \Sigma_{k}: \sigma'\neq \sigma,\\ \exists\tau\in\Sigma_{k+1}\text{ s.t. }\sigma,\sigma'\text{ are facets of } \tau }} \sgn([\sigma],\partial [\tau])\sgn([\sigma'],\partial [\tau])f([\sigma']). \end{equation}
%We now apply the constructions of Section \ref{signed}. 
We recall that a $\tau$ if it exists, is unique. \\
Our results will be obtained for this operator, and they only partially
generalize to a general $L_{k}^{up}$.

We first need to define the Cheeger constant on a simplicial complex. In fact, we shall propose four different definitions.
  \begin{enumerate}
  \item[\bf{1.}] {\bf Take multiplicities into account.}\\
   In a graph, the weight of a vertex is its degree, the number of incident edges. An edge always has weight 2, the number of its vertices. In higher dimensions, this is different.\\
   A \emph{ multiset} is a  pair $(S,m)$, where $S$ is the underlying set of the
  multiset, formed from its distinct elements, and $m:S\to\mathbb{Z}$ is an
  integer-valued function, giving the \emph{multiplicity}, possibly  negative (orientations).  $|S|
:= 
\sum_{s\in S}|m(s)|$. \\
 Let  $S$  thus be a multiset on   $\Sigma_k$ with
multiplicities in $\{-M,\ldots,0,\ldots,M\}$. Its coboundary  $\partial^*_{k+1}S$ is the  multiset of  $(k+1)$-simplices with a member of $S$ in the boundary. Each $\sigma\in \Sigma_{k+1}$ has multiplicity $\sum_{\tau\in
    \Sigma_k}m(\tau)\mathrm{sgn}([\tau],\partial[\sigma])$. 
 We put 
 $$\vol(S)=\sum_{\tau\in \Sigma_k}\deg_\tau |m(\tau)|.$$
 This makes our first definition of the 
  \emph{Cheeger constant} possible
\begin{equation}\label{eq:combinatorial-h(S_d)}h_1(\Sigma_k):=\min\limits_{\substack{S\subset_M \Sigma_k \\ S\neq \partial^*_{k}(T),\forall T\subset_M \Sigma_{k-1}}}\frac{|\partial^*_{k+1} S|}{\min\limits_{S'\ne\emptyset:\partial^*_{k+1}S'=\partial^*_{k+1}S}\vol(S')}
\end{equation}
for  sufficiently large $M$ (which then is independent of $M$, because all relations are rational).

  \item[\bf{2.}] {\bf 
    $\mathbb{Z}$-expander:}
  $$h_2(\Sigma_k):=\min\limits_{\phi\in
    C^k(\Sigma,\mathbb{Z})\setminus\mathrm{Im\,}\delta}\frac{\|\delta\phi\|_1}{\min\limits_{\psi\in\mathrm{Im\,}\delta}\|\phi+\psi\|_{1,\deg}}.$$
 %with  $\|\phi \|_{1,\deg}:=\sum_{\tau\in\Sigma_k}\deg_\tau |\phi(\tau)|$.
 Here,   $\|\phi \|_{1,\deg}:=\sum_{\tau\in
  \Sigma_k}\deg _\tau  |\phi(\tau)|$ is  the  (weighted) $l^1$-norm.

  \item[\bf{3.}] For this version of the Cheeger constant, we need the {\bf 1-Laplacian} introduced in Section \ref{plap}:
    $$h_3(\Sigma_k):=\lambda_{\min >0}(\Delta^{up}_{k,1})$$
    (smallest positive eigenvalue)

  \item[\bf{4.}] The last version of the Cheeger constant uses the {\bf filling radius:}\\
  $\|\cdot\|_{1,\deg}$ on $C^{k}(\Sigma)$ induces a quotient norm $\|\cdot\|$ on
$C^{k}(\Sigma)/\mathrm{image}(\delta_{k-1})$.  For an equivalence class $[\mathbf{x}]\in C^{k}(\Sigma)/\mathrm{image}(\delta_{k-1})$,  let $\| [\mathbf{x}]\|=\inf\limits_{x'\in [x]}\|\mathbf{x}'\|_{1,\deg}$. Then 
\begin{eqnarray*}
  h_4(\Sigma_k)&=&
                   \min\limits_{0\ne [\mathbf{x}]\in C^{k}(\Sigma,\mathbb{Z})/\mathrm{image}(\delta_{k-1})}\frac{\|\delta_k\mathbf{x}\|_1}{\| [\mathbf{x}]\|}.
\end{eqnarray*}

If ${H}^{k}(\Sigma,\R)= 0$, 
$$h_4(\Sigma_k)=\min\limits_{\mathbf{y}\in \mathrm{image}(\delta_{k})}\frac{\|\mathbf{y}\|_1}{\|\mathbf{y}\|_{\mathrm{fil}}}=\frac{1}{\max\limits_{\mathbf{y}\in \mathrm{image}(\delta_{k})}\|\mathbf{y}\|_{\mathrm{fil}}/\|\mathbf{y}\|_1}=\frac{1}{\|\delta_k^{-1}\|_{\mathrm{fil}}}$$
where $\|\mathbf{y}\|_{\mathrm{fil}}:=
\inf\limits_{x\in\delta_k^{-1} (\mathbf{y})}\|\mathbf{x}\|_{1,\deg}$ is the filling norm of $\mathbf{y}$, and\\ 
$\|\delta_k^{-1}\|_{\mathrm{fil}}:=\max\limits_{\mathbf{y}\in \mathrm{image}(\delta_{k})}\|\mathbf{y}\|_{\mathrm{fil}}/\|\mathbf{y}\|_1$ is Gromov's  filling profile \cite[Sec. 2.3]{Gromov10}.
    \end{enumerate}

   The preceding versions of the Cheeger constant look very different  from each other, but in fact they all agree.  
    \begin{Thm}[\cite{Jost/Zhang24b}]
      \begin{equation*}
        h_1=h_2=h_3=h_4=:h
      \end{equation*}
    \end{Thm}
   
 We can now formulate the Cheeger estimate.    
    \begin{Thm}[\cite{Jost/Zhang24b}]
      Suppose that $\deg_\tau>0$, $\forall \tau\in \Sigma_k$. Then
      \begin{equation*}
 \frac{h^2(\Sigma_k)}{|\Sigma_{k+1}|}\le \lambda_{\min >0}(\Delta_k^{up})\le \vol(\Sigma_k)h(\Sigma_k).      
      \end{equation*}
    \end{Thm}
Proofs of these results can also be found in \cite{Jost/Mulas/Zhang26}.

\subsubsection{The largest eigenvalue}\label{largest}
In order to control the largest eigenvalue of  $\Delta_k^{up}$, we convert a Cheeger type problem for
higher dimensional simplices into one for signed graphs.

 We   modify  the construction of the signed graph  in Section
 \ref{smallest} and use  the sign function
 \begin{equation}\label{sign1}
   s([\sigma],[\sigma'])=
   \sgn([\sigma],\partial [\tau])\sgn([\sigma'],\partial [\tau])
 \end{equation}
 for the signed graph
 $(\Gamma_k,s)$ on the vertex set $\Sigma_k$, again with $\sigma \sim \sigma'$ if they both are facets of the same (unique) $\tau \in \Sigma_{k+1}$. 
 This is the negative of the sign function from \eqref{asc2}. 
  
 The definition of the signed bipartiteness ratio \eqref{bip} now motivates a corresponding quantity.  For disjoint $A,A'\subset \Sigma_k$, let $|E^+(A,A')|=\#\{\{\sigma,\sigma'\}:\sigma\in A,\sigma'\in A',s([\sigma],[\sigma'])=1\}$  and $|E^-(A)|=\#\{\{\sigma,\sigma'\}:\sigma,\sigma'\in A,s([\sigma],[\sigma'])=-1\}$. We then define 
\begin{equation}\label{bip2}
    \beta(A,A')=\frac{2\left(|E^-(A)|+|E^-(A')|+|E^+(A,A')|\right)+|\partial(A\sqcup A')|}{\vol(A\sqcup A')}
\end{equation}
where $|\partial A|$ is the number of the edges of $(\Gamma_k,s)$ between $A$ and $\Sigma_k\setminus A$,   $\vol(A)=\sum_{\sigma\in A}\deg \sigma$ and $\deg \sigma=\#\{\tau\in \Sigma_{k+1}:\sigma \in \partial\tau\}$. 

The $j$-th  Cheeger constant   on $\Sigma_k$ is then defined as 
$$h_j(\Sigma_k)=\min\limits_{\text{disjoint } A_1,A_2,\ldots,A_{2j-1},A_{2j}\text{ in }\Sigma_k}\max\limits_{1\le i\le j}\beta(A_{2i-1},A_{2i}).$$ 
 We note that $h_j(\Sigma_k)=0$ iff $(\Gamma_k,s)$ has precisely $j$ balanced components.

This generalizes the $j$-way Cheeger
 constant $h_j(\Sigma_0)$
 of a graph  \cite{Lee12}.

In terms of this quantity, \cite{Jost/Zhang24b} shows: 
\begin{theorem}\label{thm:anti-signed-Cheeger} Let $\Sigma$  be a finite simplicial complex, and let  
$k\ge 0$. Then the largest eigenvalue $\lambda_M$ of $\Delta^{up}_k$ satisfies 
\begin{equation}\label{bip4}
\frac{ h_1(\Sigma_k)^2}{2(k+1)}\le k+2-\lambda_M(\Delta^{up}_k)\le 2h_1(\Sigma_k).
\end{equation}
 Moreover,  for any $j\ge 1$, 
\begin{equation}\label{bip5}
\frac{ h_j(\Sigma_k)^2}{Cj^6(k+1)}\le k+2-\lambda_{M+1-j}(\Delta^{up}_k)\le 2h_j(\Sigma_k)
\end{equation}
with an absolute constant $C$.
\end{theorem}

By \eqref{bip4}, $\lambda_M(\Delta^{up}_k)=k+2$ iff  $h_1(\Sigma_k)=0$. And that is the case  iff the associated  signed graph
${(\Gamma_d,s)}$ has a  balanced component, using an observation in Section \ref{sign} and observing that the sign used there 
is the opposite of the one here (so that what is antibalanced there is balanced here). 
\\

We had observed in Lemma \ref{dis} that a simplicial complex is disorientable iff the normalized up-Laplacian achieves the largest possible eigenvalue. 
As shown in \cite{Eidi24a,Eidi25}, if the $n$-complex contains no twisting cycles, then the fewer (non-branching) simple odd cycles its down dual graph possesses, the closer its largest eigenvalue of $\Delta_{n-1}^{up}$ is to $n+1$. Moreover, any $n$-simplicial complex can be made disorientable, that is, its Laplacian can attain $n+1$, through a finite number of splittings of (some of its) $n$-simplexes into two parts and the fewer such cycles exist, the fewer simplex splittings are needed to make the complex disorientable.

\section{Conclusion and Outlook}

We have connected  geometric, discrete and combinatorial structures, more precisely, Riemannian manifolds, graphs and simplicial complexes, via theories that have geometric, discrete and combinatorial versions, namely Hodge and Eckmann theory, Morse theory with its profound extensions by Witten and Floer and its discretization by Forman, the spectral theory of Laplace operators and Cheeger type inequalities and cuts. This provided us also with some opportunities to explore both structural analogies and convergence relations (when a discrete structure approximates a Riemannian manifold)  between these theories. While so far, ideas have mainly flowed from the geometric to the discrete or combinatorial setting, there are now also many prospects for the reverse direction. For instance, insight about combinatorial higher order Cheeger inequalities could pave the way for an understanding of the geometric significance of eigenvalues of higher order Hodge Laplacians. Other open problems include defining topologically informative Laplacian-based random walks on simplicial complexes that their limiting behavior could be easily analyzed, or  extending Morse–Floer-type constructions to  broader classes of dynamical systems.

The continuous versions represent some cornerstones of modern geometry. The discrete versions gain ever more relevance in data science and mathematical machine learning. Therefore, insights from the line of research described in this paper should also lead to profound advances in those fields.

This survey has traced a unified storyline across objects that are often treated in isolation: discrete and smooth Laplacians, Morse-Witten–Floer/Forman homology theories, and smooth and discrete Cheeger-type inequalities. A key message is that an integrated viewpoint is not merely philosophical, it sharpens intuition, reveals structural analogies, and enables conceptual transfer across domains.

Ultimately, viewing geometry, topology, and dynamics as different expressions of the same underlying principles offers a pathway to both deeper mathematical understanding and more principled data-driven models.\\

The entire subject feels the deep and penetrating influence of the work of Shing-Tung Yau, who advanced  the geometric side, for instance through his work on eigenvalues with Peter Li, utilized the analogies between the geometric and the discrete side in his work with Fan Chung and others, and developed a new discrete homology theory, called path homology, in a collaboration with Alexander Grigor'yan,  Yu  Muranov and Yong Lin, to name just some selected achievements.

\bibliographystyle{plain}  
\bibliography{morse.bib}

@article{Alon85,
  title={$\lambda_1$, isoperimetric inequalities for graphs, and superconcentrators},
  author={Alon, N. and Milman, V.},
  journal={Journal of Combinatorial Theory, Series B},
  volume={38},
  number={1},
  pages={73--88},
  year={1985},
  publisher={Elsevier}
}

@article{Amghibech03,
  title={{Eigenvalues of the discrete $p$-Laplacian for graphs}},
  author={Amghibech, S. },
  journal={Ars Combinatoria},
  volume={67},
  pages={283–-302},
  year={2003}
}

@article{Atay20,
  title={Cheeger constants, structural balance, and spectral clustering analysis for signed graphs},
  author={Atay, Fatihcan  and Liu, Shiping},
  journal={Discrete Mathematics},
  volume={343},
  number={1},
  pages={111616},
  year={2020}
}

@article{Bannwart24,
  title={About non-uniqueness when removing closed orbits in {Morse-Smale} vector fields},
  author={Bannwart, Clemens},
  journal={arXiv preprint arXiv:2410.02363},
  year={2024}
}

@book {Banyaga04,

    AUTHOR = {Augustin Banyaga and David Hurtubise},

     TITLE = {Lectures on {M}orse homology},

 PUBLISHER = {Kluwer},

      YEAR = {2004},

}

@article {Bauer13,

    AUTHOR = {Frank Bauer and J\"urgen Jost},

     TITLE = {Bipartite and neighborhood graphs and the spectrum of the normalized graph Laplacian},

   JOURNAL = {Comm. Anal. Geom.},

    VOLUME = {21},

      YEAR = {2013},

      PAGES = {787 -- 845},

}

@ARTICLE{Belkin03,

  author={Belkin, Mikhail and Niyogi, Partha},

  journal={Neural Computation}, 

  title={Laplacian Eigenmaps for Dimensionality Reduction and Data Representation}, 

  year={2003},

  volume={15},

  number={6},

  pages={1373--1396},

  doi={10.1162/089976603321780317}}

@article{Belkin08,
title = {Towards a theoretical foundation for {L}aplacian-based manifold methods},
journal = {Journal of Computer and System Sciences},
volume = {74},
number = {8},
pages = {1289-1308},
year = {2008},
note = {Learning Theory 2005},
issn = {0022-0000},
doi = {https://doi.org/10.1016/j.jcss.2007.08.006},
url = {https://www.sciencedirect.com/science/article/pii/S0022000007001274},
author = {Mikhail Belkin and Partha Niyogi}
}

@article{Beltrami64,
  title={Ricerche di analisi applicata alla geometria},
  author={Beltrami, Eugenio},
  journal={Giornale di Matematiche},
  volume={},
  number={},
  pages={vol. II: 267--282, 297--306, 331--339, 355--?; vol. III: 15--22, 33--45, 82--91, 228--240, 311--314},
  year={1864}
}

@book{Beltrami02,
  title={{Opere matematiche di Eugenio Beltrami. Pubblicate per cura della Facolt\`a di scienze della R. Universit\`a di Roma}},
  author={Beltrami, Eugenio},
  year={1902--20},
  publisher={U. Hoepli},
address={Milano}
}

@ARTICLE{Benedetti10,
       author = {{Benedetti}, Bruno},
        title = "{Discrete Morse theory is at least as perfect as Morse theory}",
      journal = {arXiv e-prints},
         year = 2010,
        month = oct,
          eid = {arXiv:1010.0548},
        pages = {arXiv:1010.0548},
archivePrefix = {arXiv},
       eprint = {1010.0548},
 primaryClass = {math.DG},
       adsurl = {https://ui.adsabs.harvard.edu/abs/2010arXiv1010.0548B}
}

@article{Benedetti16,
  title={Smoothing discrete {M}orse theory},
  author={Benedetti, Bruno},
  journal={Ann. Sc. Norm. Super. Pisa Cl. Sci. (5)},
  pages={335--368},
  year={2016}
}

@article{Brehm87,
  title={Combinatorial manifolds with few vertices},
  author={Brehm, Ulrich and K{\"u}hnel, Wolfgang},
  journal={Topology},
  volume={26},
  number={4},
  pages={465--473},
  year={1987}
}

@article{Buser78,
  title={{\"Uber eine Ungleichung von Cheeger}},
  author={Buser, Peter},
  journal={Mathematische Zeitschrift},
  volume={158},
  number={3},
  pages={245--252},
  year={1978}
}

@article {Buser80,

    AUTHOR = {Peter Buser},

     TITLE = {{On Cheeger's inequality: $\lambda_1 \ge \frac{1}{4}h^2$}},

   JOURNAL = {AMS Proc.Symp.Pure Math.},

    VOLUME = {36},

      YEAR = {1980},

    NUMBER = {},

     PAGES = {29--77},

}

@article {Buser82,
    AUTHOR = {Buser, Peter},
     TITLE = {A note on the isoperimetric constant},
   JOURNAL = {Ann. Sci. \'Ecole Norm. Sup. (4)},
  FJOURNAL = {Annales Scientifiques de l'\'Ecole Normale Sup\'erieure.
              Quatri\`eme S\'erie},
    VOLUME = {15},
      YEAR = {1982},
    NUMBER = {2},
     PAGES = {213--230},
      ISSN = {0012-9593},
   MRCLASS = {58G25 (52A40 53C20)},
  MRNUMBER = {683635},
MRREVIEWER = {Scott\ Wolpert},
       URL = {http://www.numdam.org/item?id=ASENS_1982_4_15_2_213_0},
}

@incollection {Cheeger70,

    AUTHOR = {Jeff Cheeger},

     TITLE = {A lower bound for the smallest eigenvalue of the {L}aplacian},

 BOOKTITLE = {Problems in Analysis},

     PAGES = {195--199},

 PUBLISHER = {Princeton Univ.Press},

   ADDRESS = {},

    EDITOR = {},

      YEAR = {1970},

}

@book {Chung97,
    AUTHOR = {Chung, Fan R. K.},
     TITLE = {Spectral graph theory},
    SERIES = {CBMS Regional Conference Series in Mathematics},
    VOLUME = {92},
 PUBLISHER = {Conference Board of the Mathematical Sciences, Washington, DC;
              by the American Mathematical Society, Providence, RI},
      YEAR = {1997},
     PAGES = {xii+207},
      ISBN = {0-8218-0315-8},
   MRCLASS = {58G99 (05C50 35P05 46N20 47N20)},
  MRNUMBER = {1421568},
MRREVIEWER = {Robert\ Brooks},
}

@book {Conley78,

    AUTHOR = {Charles Conley},

     TITLE = {{Isolated invariant sets and the Morse index}},

 PUBLISHER = {CBMS Reg. Conf.
Ser. Math. 38, AMS},

   ADDRESS = {Providence, R.I.},

      YEAR = {1978},

}

@article{DeGiorgi75,
  title={Su un tipo di convergenza variazionale},
  author={Ennio De Giorgi and Tullio Franzoni},
  journal={Atti
della Accademia Nazionale dei Lincei. Classe di Scienze Fisiche, Matematiche e Naturali},
  volume={58},
  number={6},
  pages={842--850},
  year={1975}
}

@article{Degiovanni10,
  title={On topological and metric critical point theory},
  author={Degiovanni, Marco},
  journal={Journal of Fixed Point Theory and Applications},
  volume={7},
  number={1},
  pages={85--102},
  year={2010}
}

@article{Degiovanni94,
  title={A critical point theory for nonsmooth functional},
  author={Degiovanni, Marco and Marzocchi, Marco},
  journal={Annali di Matematica Pura ed Applicata},
  volume={167},
  number={1},
  pages={73--100},
  year={1994}
}

@article {Deidda23,
    AUTHOR = {Deidda, Piero and Putti, Mario and Tudisco, Francesco},
     TITLE = {Nodal domain count for the generalized graph
              {$p$}-{L}aplacian},
   JOURNAL = {Appl. Comput. Harmon. Anal.},
  FJOURNAL = {Applied and Computational Harmonic Analysis. Time-Frequency
              and Time-Scale Analysis, Wavelets, Numerical Algorithms, and
              Applications},
    VOLUME = {64},
      YEAR = {2023},
     PAGES = {1--32},
      ISSN = {1063-5203,1096-603X},
   MRCLASS = {47J10 (35J10)},
  MRNUMBER = {4530637},
MRREVIEWER = {Abdallah\ Maichine},
       DOI = {10.1016/j.acha.2022.12.003},
       URL = {https://doi.org/10.1016/j.acha.2022.12.003},
}

@article {DodziukPatodi76,
    AUTHOR = {Dodziuk, Jozef and Patodi, V. K.},
     TITLE = {Riemannian structures and triangulations of manifolds},
   JOURNAL = {J. Indian Math. Soc. (N.S.)},
  FJOURNAL = {The Journal of the Indian Mathematical Society. New Series},
    VOLUME = {40},
      YEAR = {1976},
    NUMBER = {1-4},
     PAGES = {1--52},
      ISSN = {0019-5839,2455-6475},
   MRCLASS = {58G10 (57D20)},
  MRNUMBER = {488179},
MRREVIEWER = {M.\ L.\ Gromov},
}

@article {Dodziuk76,
    AUTHOR = {Dodziuk, Jozef},
     TITLE = {Finite-difference approach to the {H}odge theory of harmonic
              forms},
   JOURNAL = {Amer. J. Math.},
  FJOURNAL = {American Journal of Mathematics},
    VOLUME = {98},
      YEAR = {1976},
    NUMBER = {1},
     PAGES = {79--104},
      ISSN = {0002-9327,1080-6377},
   MRCLASS = {58A10 (58G99)},
  MRNUMBER = {407872},
MRREVIEWER = {D.\ B.\ Fuchs},
       DOI = {10.2307/2373615},
       URL = {https://doi.org/10.2307/2373615},
}

@article{Dodziuk84,
  title={Difference equations, isoperimetric inequality and transience of certain random walks},
  author={Dodziuk, Jozef},
  journal={Transactions of the American Mathematical Society},
  volume={284},
  number={2},
  pages={787--794},
  year={1984}
}

@article{Dotterrer12,
  title={Coboundary expanders},
  author={Dotterrer, Dominic and Kahle, Matthew},
  journal={Journal of Topology and Analysis},
  volume={4},
  number={04},
  pages={499--514},
  year={2012}
}

@article{Eckmann44,
author = {Eckmann, Beno},
issn = {0010-2571},
journal = { Comment. Math. Helv.},
month = dec,
number = {1},
pages = {240--255},
title = {{Harmonische Funktionen und Randwertaufgaben in einem Komplex}},
volume = {17},
year = {1944}
}

@book{Edelsbrunner10,
  title={Computational topology: an introduction},
  author={Edelsbrunner, Herbert and Harer, John},
  year={2010},
  publisher={American Mathematical Soc.}
}

@article{Eidi21,
  title={{Floer homology: From generalized Morse-Smale dynamical systems to Forman's combinatorial vector fields}},
  author={Eidi, Marzieh and Jost, J{\"u}rgen},
  journal={arXiv preprint arXiv:2105.02567},
  year={2021}
}

@article{Eidi23,
  title={Irreducibility of {M}arkov Chains on simplicial complexes, the Spectrum of the Discrete {Hodge Laplacian} and Homology},
  author={Eidi, Marzieh and Mukherjee, Sayan},
  journal={arXiv preprint arXiv:2310.07912},
  year={2023}
}

@article{Eidi24,
author = {Marzieh Eidi and J\"urgen Jost},
title = {{Floer homology\,:\, from generalized Morse-Smale dynamical systems to Forman's combinatorial
vector fields}},
doi = {10.1007/s40304-022-00314-6},
journal = {Communications in Mathematics and Statistics},
pages = {695--720},
year = {2024},
volume = {12},
number = {4},
issn = {2194-6701},
}

@article{Eidi24a,
  title={Higher Order Bipartiteness vs Bi-Partitioning in Simplicial Complexes},
  author={Eidi, Marzieh and Mukherjee, Sayan},
  journal={arXiv preprint arXiv:2409.00682},
  year={2024}
}

@InProceedings{Eidi25,
  author =	{Eidi, Marzieh and Mukherjee, Sayan},
  title =	{{Higher Order Bipartiteness vs Bi-Partitioning in Simplicial Complexes}},
  booktitle =	{41st International Symposium on Computational Geometry (SoCG 2025)},
  pages =	{45:1--45:12},
  series =	{Leibniz International Proceedings in Informatics (LIPIcs)},
  ISBN =	{978-3-95977-370-6},
  ISSN =	{1868-8969},
  year =	{2025},
  volume =	{332},
  editor =	{Aichholzer, Oswin and Wang, Haitao},
  publisher =	{Schloss Dagstuhl -- Leibniz-Zentrum f{\"u}r Informatik},
  address =	{Dagstuhl, Germany},
  URL =		{https://drops.dagstuhl.de/entities/document/10.4230/LIPIcs.SoCG.2025.45},
  URN =		{urn:nbn:de:0030-drops-231972},
  doi =		{10.4230/LIPIcs.SoCG.2025.45}
}

@article {Floer88c,

    AUTHOR = {Andreas Floer},

     TITLE = {A relative {M}orse index for the symplectic action},

   JOURNAL = {Comm.Pure Appl.Math.},

    VOLUME = {41},

      YEAR = {1988},

    NUMBER = {},

     PAGES = {393--407},

}

@article {Floer89,

    AUTHOR = {Andreas Floer},

     TITLE = {Witten's complex and infinite dimensional {M}orse theory},

   JOURNAL = {{\it J. Diff. Geom.}},

    VOLUME = {30},

      YEAR = {1989},

    NUMBER = {},

     PAGES = {207--221},

}

@article{Forman98a,
title = {Morse Theory for Cell Complexes},
journal = {Advances in Mathematics},
volume = {134},
number = {1},
pages = {90-145},
year = {1998},
issn = {0001-8708},
doi = {https://doi.org/10.1006/aima.1997.1650},
url = {https://www.sciencedirect.com/science/article/pii/S0001870897916509},
author = {Robin Forman}
}

@article{Forman98b,
  title={Combinatorial vector fields and dynamical systems},
  author={Robin Forman},
  journal={Mathematische Zeitschrift},
  year={1998},
  volume={228},
  pages={629-681},
}

@article{Forman98c,
  title={{Witten-Morse theory for cell complexes}},
  author={Robin Forman},
  journal={Topology},
  year={1998},
  volume={37},
  pages={945--979},
}

@article{Franks79,
  title={Morse-Smale flows and homotopy theory},
  author={J. Franks},
  journal={Topology},
  year={1979},
  volume={18},
  pages={199-215}
}

@article{Gallais10,
  title={{Combinatorial realization of the Thom-Smale complex via discrete Morse theory}},
  author={Gallais, {\'E}tienne},
  journal={Annali della Scuola Normale Superiore di Pisa-Classe di Scienze},
  volume={9},
  number={2},
  pages={229--252},
  year={2010}
}

@article{Garcia16,
  title={Continuum limit of total variation on point clouds},
  author={Garc{\'\i}a Trillos, Nicol{\'a}s and Slep{\v{c}}ev, Dejan},
  journal={Archive for Rational Mechanics and Analysis},
  volume={220},
  number={1},
  pages={193--241},
  year={2016}
}

@article{Garcia20,
  title={{Error estimates for spectral convergence of the graph Laplacian on random geometric graphs toward the Laplace--Beltrami operator}},
  author={Garc{\'\i}a Trillos, Nicol{\'a}s and Gerlach, Moritz and Hein, Matthias and Slep{\v{c}}ev, Dejan},
  journal={Foundations of Computational Mathematics},
  volume={20},
  number={4},
  pages={827--887},
  year={2020}
}

@article{Grigoryan14,
  title={Graphs associated with simplicial complexes},
  author={Grigor'yan, Alexander and Muranov, Yu V and Yau, Shing-Tung},
journal={Homology, Homotopy and Applications},
  volume={16},
  number={1},
  pages={295--311},
  year={2014}
}

@article{Gromov10,
  title={{Singularities, expanders and topology of maps. Part 2: From combinatorics to topology via algebraic isoperimetry}},
  author={Gromov, Mikhail},
  journal={Geometric and Functional Analysis},
  volume={20},
  number={2},
  pages={416--526},
  year={2010}
}

@Article{Harary55,
  author = 	 {Harary, F.},
  title = 	 {On the notion of balance of a signed graph},
  journal = 	 {Michigan Math. J.},
  year = 	 {1955},
  OPTkey = 	 {},
  volume = 	 {2},
  OPTnumber = 	 {},
  pages = 	 {143--146},
  OPTmonth = 	 {},
  OPTnote = 	 {},
  OPTannote = 	 {}
}

@Article{Harary57,
  author = 	 {Harary, F.},
  title = 	 {Structural duality},
  journal = 	 {Behavioral Sci.},
  year = 	 {195},
  OPTkey = 	 {},
  volume = 	 {2},
  OPTnumber = 	 {},
  pages = 	 {255--265},
  OPTmonth = 	 {},
  OPTnote = 	 {},
  OPTannote = 	 {}
}

@Article{Harary59,
  author = 	 {Harary, F.},
  title = 	 {On the measurement of structural balance},
  journal = 	 {Behavioral Sci.},
  year = 	 {1959},
  OPTkey = 	 {},
  volume = 	 {4},
  OPTnumber = 	 {},
  pages = 	 { 316--323},
  OPTmonth = 	 {},
  OPTnote = 	 {},
  OPTannote = 	 {}
}

@article {Harary60,
    AUTHOR = {Harary, F. and Norman, R. Z.},
     TITLE = {Some properties of line digraphs},
   JOURNAL = {Rend. Circ. Mat. Palermo (2)},
  FJOURNAL = {Rendiconti del Circolo Matematico di Palermo. Serie II},
    VOLUME = {9},
      YEAR = {1960},
     PAGES = {161--168},
    }

@ARTICLE {Hein10,
    author  = "Hein, M. and {B\"uhler}, T.",
    title   = "{An inverse power method for nonlinear eigenproblems with applications in 1--spectral clustering and sparse {PCA}}",
    journal = "NIPS",
    year    = "2010",
    pages   = "847--855"
}

@article{Hilbert04,
  title={{Grundz{\"u}ge einer allgemeinen Theorie der linearen Integralgleichungen. (Erste Mitteilung)}},
  author={Hilbert, David},
  journal={Nachrichten  der Gesellschaft der Wissenschaften zu G{\"o}ttingen, Mathematisch-Physikalische Klasse},
  volume={1904},
  pages={49--91},
  year={1904}
}

@book {Hodge41,

    AUTHOR = {William Hodge},

     TITLE = {The Theory and Applications of Harmonic Integrals},

 PUBLISHER = {Cambr.Univ.Press},

   ADDRESS = {},

      YEAR = {1941, 2nd ed. 1952},

}

@Article{Horak13,
  author = 	 {Danijela Horak and  J\"urgen Jost},
  title = 	 {Spectra of combinatorial {L}aplace operators on simplicial
complexes},
  journal = 	 {Adv. Math.},
  year = 	 {2013},
  OPTkey = 	 {},
  volume = 	 {244},
  OPTnumber = 	 {},
  pages = 	 {303--336},
  OPTmonth = 	 {},
  OPTnote = 	 {},
  OPTannote = 	 {}
}

@article{Ioffe96,
  title={Metric critical point theory 1. Morse regularity and homotopic stability of a minimum},
  author={Ioffe, Alexander and Schwartzman, Efim},
  journal={Journal de mathematiques pures et appliqu{\'e}es},
  volume={75},
  number={2},
  pages={125--154},
  year={1996}
}

@book{Joharinad23,

    AUTHOR = {Joharinad, Parvaneh and Jost, J\"urgen},

     TITLE = {Mathematical Principles of Topological and Geometric Data Analysis},
SERIES ={Mathematics of Data},
 PUBLISHER = {Springer},

   ADDRESS = {},

      YEAR = {2023},

}

@book{Jost14,
  title={Mathematical methods in biology and neurobiology},
  author={Jost, J{\"u}rgen},
  year={2014},
  publisher={Springer}
}

@book {Jost17,

    AUTHOR = {J\"urgen Jost},

     TITLE = {{Riemannian geometry and geometric analysis}},

 PUBLISHER = {Springer},

   ADDRESS = {},

      YEAR = {8th ed., 2026},

}

@book {Jost25,

    AUTHOR = {Jost, J\"urgen},

     TITLE = {{Bernhard Riemann. On the hypotheses which lie at the bases of geometry}},

 PUBLISHER = {Birkh\"auser},
 

   SERIES = {Classic Texts in the Sciences},

      YEAR = {2nd ed. 2025}

}

@book{Jost98,
  title={Calculus of variations},
  author={Jost, J{\"u}rgen and Li-Jost, Xianqing},
  volume={64},
  year={1998},
  publisher={Cambridge University Press}
}

@book{Jost/Mulas/Zhang26,

    AUTHOR = {Jost, J\"urgen and Mulas, Raffaella and Zhang, Dong},

     TITLE = {Spectra of discrete structures},
SERIES ={Cambridge Studies in Advanced Mathematics},
 PUBLISHER = {Cambridge University Press},

   ADDRESS = {},

      YEAR = {2026},

}

@article{Jost/Zhang21a,
  title={Discrete-to-Continuous Extensions: piecewise multilinear extension, min-max theory and spectral theory},
  author={Jost, J{\"u}rgen and Zhang, Dong},
  journal={arXiv preprint arXiv:2106.04116},
  year={2021}
}

@article{Jost/Zhang21b,
  title={Discrete-to-Continuous Extensions: Lov$\backslash$'asz extension, optimizations and eigenvalue problems},
  author={Jost, J{\"u}rgen and Zhang, Dong},
  journal={arXiv preprint arXiv:2106.03189},
  year={2021}
}

@article{Jost/Zhang24a,
  title={Discrete-to-Continuous Extensions: {L}ov{\'a}sz Extension and {M}orse Theory},
  author={Jost, J{\"u}rgen and Zhang, Dong},
  journal={Discrete \& Computational Geometry},
  volume={72},
  number={1},
  pages={49--72},
  year={2024}
}

@article{Jost/Zhang24b,
  title={Cheeger inequalities on simplicial complexes},
  author =	 {Jost, J{\"u}rgen and Zhang, Dong},
  journal =	 {Ann. Sc. Norm. Super. Pisa Cl. Sci. (5)},
doi = {https://doi.org/10.2422/2036-2145.202307_009}, URL = {https://doi.org/10.2422/2036-2145.202307_009},
year = {2024}
}

@article{Katriel94,
  title={Mountain pass theorems and global homeomorphism theorems},
  author={Katriel, Guy},
  journal={Annales de l'Institut Henri Poincar{\'e} C, Analyse non lin{\'e}aire},
  volume={11},
  number={2},
  pages={189--209},
  year={1994}
}

@incollection{Kuhnel90,
  title={Triangulations of manifolds with few vertices},
  author={K{\"u}hnel, Wolfgang},
  booktitle={Advances in differential geometry and topology},
  pages={59--114},
  year={1990},
  publisher={World Scientific}
}

@incollection {Lee12,
    AUTHOR = {Lee, J. R. and Oveis Gharan, S. and Trevisan, L.},
     TITLE = {Multi-way spectral partitioning and higher-order {C}heeger
              inequalities},
 BOOKTITLE = {S{TOC}'12---{P}roceedings of the 2012 {ACM} {S}ymposium on
              {T}heory of {C}omputing},
     PAGES = {1117--1130},
 PUBLISHER = {ACM, New York},
      YEAR = {2012},
   MRCLASS = {05C40 (05C50 05C85)},
  MRNUMBER = {2961569},
       DOI = {10.1145/2213977.2214078},
       URL = {http://dx.doi.org/10.1145/2213977.2214078},
}

@article {Li80a,

    AUTHOR = {Peter Li and Shing-Tung Yau},

     TITLE = {{Estimates of eigenvalues of a
compact Riemannian manifold}},

   JOURNAL = {AMS Proc.Symp.Pure Math.},

    VOLUME = {36},

      YEAR = {1980},

    NUMBER = {},

     PAGES = {205--240},
}

@article{Li83,
  title={On the {S}chr{\"o}dinger equation and the eigenvalue problem},
  author={Li, Peter and Yau, Shing-Tung},
  journal={Communications in Mathematical Physics},
  volume={88},
  number={3},
  pages={309--318},
  year={1983}
}

@article {Lickorish,
    AUTHOR = {Lickorish, W. B. R.},
     TITLE = {Unshellable triangulations of spheres},
   JOURNAL = {European J. Combin.},
  FJOURNAL = {European Journal of Combinatorics},
    VOLUME = {12},
      YEAR = {1991},
    NUMBER = {6},
     PAGES = {527--530},
      ISSN = {0195-6698,1095-9971},
   MRCLASS = {57Q15 (52B70 57Q45)},
  MRNUMBER = {1136394},
MRREVIEWER = {C.\ Kearton},
       DOI = {10.1016/S0195-6698(13)80103-5},
       URL = {https://doi.org/10.1016/S0195-6698(13)80103-5},
}

@article{Morse25,
 URL = {http://www.jstor.org/stable/1989110},
 author = {Marston Morse},
 journal = {Transactions of the American Mathematical Society},
 number = {3},
 pages = {345--396},
 title = {Relations Between the Critical Points of a Real Function of $n$ Independent Variables},
 volume = {27},
 year = {1925}
}

@ARTICLE{Polya51,
  author = 	 {P\'olya, G. and Szeg\"o, S.},
  title = 	 "{Isoperimetric inequalities in mathematical physics}",
  journal = {Annals of Math. Studies},
  year = {1951},
  volume = {27}
  }

@article{Rosenthal17,
author = {Parzanchevski, Ori and Rosenthal, Ron},
title = {Simplicial complexes: Spectrum, homology and random walks},
journal = {Random Structures \& Algorithms},
volume = {50},
number = {2},
pages = {225-261},
doi = {https://doi.org/10.1002/rsa.20657},
url = {https://onlinelibrary.wiley.com/doi/abs/10.1002/rsa.20657},
eprint = {https://onlinelibrary.wiley.com/doi/pdf/10.1002/rsa.20657},
year = {2017},
}

@book {Schwarz93,

    AUTHOR = {Matthias Schwarz},

     TITLE = {Morse homology},

 PUBLISHER = {Birkh\"auser},

   ADDRESS = {},

      YEAR = {1993},

}

@article{Steenbergen14,
  title        = {A Cheeger-Type Inequality on Simplicial Complexes},
  author       = {Steenbergen, John and Klivans, Caroline and Mukherjee, Sayan},
  journal      = {Advances in Applied Mathematics},
  volume       = {56},
  pages        = {56--77},
  year         = {2014},
  doi          = {10.1016/j.aam.2014.01.002},
  url          = {https://doi.org/10.1016/j.aam.2014.01.002},
  issn         = {0196-8858}
}

@article{Trevisan12,
  title={Max cut and the smallest eigenvalue},
  author={Trevisan, L.},
  journal={SIAM J. Comput.},
  volume={41},
  number={6},
  pages={1769--1786},
  year={2012}
}

@Article{Weyl12,
  author = 	 {Weyl, Hermann},
  title = 	 {{Das asymptotische Verteilungsgesetz der Eigenwerte linearer partieller Differentialgleichungen}},
  journal = 	 {Math.Ann.},
  year = 	 {1912},
  OPTkey = 	 {},
  volume = 	 {71},
  OPTnumber = 	 {},
  pages = 	 {441--469},
  OPTmonth = 	 {},
  OPTnote = 	 {},
  OPTannote = 	 {}
}

@article {Witten82a,

    AUTHOR = {Edward Witten},

     TITLE = {Constraints on supersymmetry breaking},

   JOURNAL = {Nucl. Phys. B},

    VOLUME = {202},

      YEAR = {1982},

    NUMBER = {},

     PAGES = {253--316},

}

@article {Witten82,

    AUTHOR = {Edward Witten},

     TITLE = {Supersymmetry and {M}orse theory},

   JOURNAL = {J. Diff. Geom},

    VOLUME = {17},

      YEAR = {1982},

    NUMBER = {},

     PAGES = {661--692},

}

@article{Yaptieu17,
   title = {{Discrete Morse-Bott theory for CW complexes}},
   author = {Yaptieu, Sylvia},
   journal = {arXiv preprint arXiv:1711.10983},
    year = {2017}
}

@article {Zhang25,
    AUTHOR = {Zhang, Dong},
     TITLE = {Homological eigenvalues of graph {$p$}-{L}aplacians},
   JOURNAL = {J. Topol. Anal.},
  FJOURNAL = {Journal of Topology and Analysis},
    VOLUME = {17},
      YEAR = {2025},
    NUMBER = {2},
     PAGES = {555--606},
      ISSN = {1793-5253,1793-7167},
   MRCLASS = {47J10 (05C50 31C20 35J92 35P05 35R02 55N31 58E05)},
  MRNUMBER = {4899943},
       DOI = {10.1142/S1793525323500346},
       URL = {https://doi.org/10.1142/S1793525323500346},
}

@article{Hein2007,
  author    = {Matthias Hein and Jean-Yves Audibert and Ulrike von Luxburg},
  title     = {Graph Laplacians and their Convergence on Random Neighborhood Graphs},
  journal   = {Journal of Machine Learning Research},
  volume    = {8},
  pages     = {1325--1368},
  year      = {2007},
  url       = {http://www.jmlr.org/papers/volume8/hein07a/hein07a.pdf}
}

@incollection{Lovasz1993,
  author    = {L{\'a}szl{\'o} Lov{\'a}sz},
  title     = {Random Walks on Graphs: A Survey},
  booktitle = {Combinatorics, Paul Erd{\H o}s is Eighty, Vol. 2},
  editor    = {D. Mikl{\'o}s and V.T. S{\'o}s and T. S{\'o}s},
  volume    = {2},
  pages     = {1--46},
  year      = {1993},
  publisher = {J{\'a}nos Bolyai Mathematical Society},
  address   = {Keszthely, Hungary},
  note      = {Dedicated to the marvelous random walk of Paul Erd{\H o}s},
  url       = {https://www.cs.yale.edu/publications/techreports/tr1029.pdf}
}

@article{Mukherjee2016,
  author    = {Sayan Mukherjee and John Steenbergen},
  title     = {Random walks on simplicial complexes and harmonics},
  journal   = {Random Structures \& Algorithms},
  volume    = {49},
  number    = {3},
  pages     = {379--405},
  year      = {2016},
  doi       = {10.1002/rsa.20645},
  url       = {https://doi.org/10.1002/rsa.20645}
}

@book{Hsu2002,
  author    = {Elton P. Hsu},
  title     = {Stochastic Analysis on Manifolds},
  series    = {Graduate Studies in Mathematics},
  volume    = {38},
  publisher = {American Mathematical Society},
  address   = {Providence, RI},
  year      = {2002},
  isbn      = {978-0-8218-3210-8},
  doi       = {10.1090/gsm/038},
  url       = {https://bookstore.ams.org/gsm-38},
  note      = {American Mathematical Society, Graduate Studies in Mathematics, Vol.~38}
}

@book{Milnor,
 ISBN = {9780691079967},
 URL = {http://www.jstor.org/stable/j.ctt183psc9},
 author = {John Milnor},
 publisher = {Princeton University Press},
 title = {Lectures on the H-Cobordism Theorem},
 urldate = {2024-11-16},
 year = {1965}
}

@misc{cancel,
      title={A proof of Morse's theorem about the cancellation of critical points}, 
      author={Francois Laudenbach},
      year={2013},
      eprint={1307.2545},
      archivePrefix={arXiv},
      primaryClass={math.GT},
      url={https://arxiv.org/abs/1307.2545}, 
}

\end{document}